\documentclass[11pt, oneside]{article}   	
\usepackage{geometry}              		
\usepackage{color}  
\usepackage{booktabs}
\usepackage{longtable}
\usepackage{rotating}
\usepackage{soul}
\usepackage[all]{xy}
\usepackage{graphicx}				
\usepackage{amssymb}
\usepackage{amsmath}
\usepackage{hyperref}

\title{
Perfectoid Diamonds and $n$-Awareness. A Meta-Model of Subjective Experience.}
\author{Shanna Dobson$^{1}$ and Robert Prentner$^{2}$}\date{}							
\begin{document}
\maketitle 
\vspace{-0.66cm}
\begin{footnotesize}
\noindent $1$ Department of Mathematics, California State University, Los Angeles, CA.\\
$2$ Center for the Future Mind, Florida Atlantic University, Boca Raton, FL.\\
\end{footnotesize}
\setstcolor{red}
\begin{center}
\subsection*{Abstract}
\end{center}
In this paper, we propose a mathematical model of subjective experience in terms of  classes of hierarchical geometries of representations (``$n$-awareness''). 

We first outline a general framework by recalling concepts from higher category theory, homotopy theory, and the theory of $(\infty,1)$-topoi. We then state three conjectures that enrich this framework. We first propose that the $(\infty,1)$-category of a geometric structure known as perfectoid diamond is an $(\infty,1)$-topos. In order to construct a topology on the $(\infty,1)$-category of diamonds we then propose that
topological localization, in the sense of Grothendieck-Rezk-Lurie $(\infty,1)$-topoi, extends to the $(\infty,1)$-category of diamonds.  We provide a small-scale model using triangulated categories.  Finally, our meta-model takes the form of Efimov K-theory of the  $(\infty,1)$-category of perfectoid diamonds, which illustrates structural equivalences between the category of diamonds and subjective experience (i.e. its privacy, self-containedness, and self-reflexivity). 

Based on this, we investigate implications of the model. We posit a grammar (“$n$-declension”) for a novel language to express  $n$-awareness, accompanied by a new  temporal scheme (``$n$-time''). Our framework allows us to revisit old problems in the philosophy of time: how is change possible and what do we mean by simultaneity and coincidence? We also examine the notion of ``self'' within our framework.
A new model of personal identity is introduced which resembles a categorical version of the ``bundle theory'': selves are not substances in which properties inhere but  (weakly) persistent moduli spaces in the K-theory of perfectoid diamonds.
%
%
\\\\
\noindent \textbf{Keywords}: Perfectoid spaces, representation, $(\infty,1)$-topoi, higher category theory, homotopy theory,  $n$-time, simultaneity, selves, localization, Efimov K-theory.\\

\section{Introduction}

\subsection{Representing Awareness as Diamond}
``Representationalism'' in the philosophy of mind is the well-known view that regards our thoughts to be representations of external things, perceiving a ``mental image'' rather than directly accessing their nature. We are most interested in representationalist theories \textit{of consciousness} \cite{Tye09}. There are many variants of such theories, and we refrain of a more detailed overview, let alone listing and comparing their relative merits and problems.  However, representationalism about consciousness comes with several well-known problems that have to do with a perceived mismatch between subjective experiences and representations:

\textit{Privacy} refers to the inaccessibility of subjective experience ``from the outside''. The ``inside'' of an experience is not directly accessible, but only reflected by the outside. A metaphor that illustrates this is that of a (mineralogical) diamond. A diamond displays its (invisible) impurities as (visible) colors. 
The following properties particularly apply to the experience of selves:
\textit{Self-containedness} refers to the property that any possible experiential pattern is ``always already'' \cite{Hei27} prefigured within the totally of an experience.\footnote{This property further carries the connotation of being a ``substance'' in the sense of an independently existing entity (e.g. \cite{Spi77}).} \textit{Self-reflexivity} refers to the idea that all parts that make up an experience are mirrored in their whole experiential context. 

We propose a concrete mathematical structure, the ``perfectoid diamond'' \cite{SW20}, as a model of subjective experience and appeal to Efimov K-theory \cite{Hoy18} to study the equivalence relations across all such diamonds. 
This would eventually lead to a geometric unification of models that study (all possible forms of) subjective experience, which is why we refer to it as our ``meta-model'', analogous to a grand unified theory of mathematical structures (the ``Langlands-program'' \cite{Far16}).

We will first outline a basic (``vanilla'') framework of $(\infty,1)$-topoi, followed by three associated conjectures that define our meta-model:

\begin{itemize}

\item \textbf{Conjecture 1}. The  $(\infty,1)$-category of perfectoid diamonds is an $(\infty,1)$-topos.  

\item \textbf{Conjecture 2}. Topological localization, in the sense of Grothendieck-Rezk-Lurie $(\infty,1)$-topoi, extends to the $(\infty, 1)$-category of diamonds.

\item \textbf{Conjecture 3}.  The meta-model takes the form of Efimov K-theory of the $(\infty,1)$-category of diamonds.
\end{itemize}

\begin{figure}[htb]
\begin{tabular}{p{0.35\textwidth} p{0.4\textwidth}}
\begin{minipage}{1\textwidth}
\includegraphics[width=5cm]{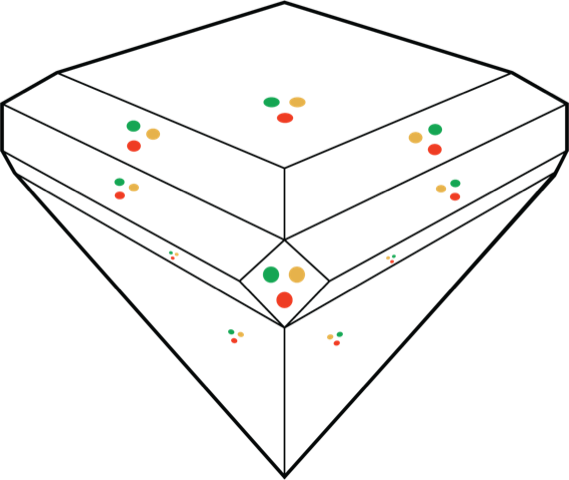}
\end{minipage}
&
\begin{tabular}{p{0.22\textwidth} p{0.26\textwidth}}
\toprule
privacy &  geometric points \& pro-\'etale topology \\
self-containedness & descent, presentabillity\\
self-reflexivity & localization via quasi-isomorphisms\\
\bottomrule
\end{tabular} 
\end{tabular}
\caption{Subjective experience as a perfectoid diamond. 
The interior points of a diamond resemble mineralogical impurities made visible as reflections on its surface. The impurity is never detected directly, but only indirectly through its surface reflections. In the figure, each impurity gives rise to three different colors (red-green-gold), depending on the ``angle'' they are viewed from and symbolizing different representations of the interior. The number is arbitrary.
The main analogies between commonly assumed properties of subjective experience and their perfectoid diamond correlation are shown on the right hand side of the figure.
}  
\label{fig:diamond}
\end{figure}

The reality of experience is neither reducible nor irreducible to its representation. Nothing dies out in the perfectoid structure.
The main proposed analogies between subjective experience and perfectoid diamonds are  stated in Fig. \ref{fig:diamond}. 

In brief, a diamond is a sheaf $Y$ on the category of perfectoid spaces \cite{per} of characteristic $p$, which is constructed as the quotient of a perfectoid space by a pro-\'etale equivalence relation. Specifically, this quotient lives in a category of sheaves on the site of perfectoid spaces with pro-\'etale covers. Perfectoid spaces are analytic spaces over nonarchimedean fields, variants of Huber space \cite{hub}. The pro-\'etale topology is a topology on the category 
of perfectoid spaces. A pro-\'etale equivalance relation is a way of gluing together perfectoid spaces under the pro-\'etale topology. A sheaf is a tool used to pack together local data on a topological space. 
Diamonds contain impurities which are profinitely many copies of geometric points%
\footnote{Let $S$ be a scheme defined over a field $k$ and equipped with a morphism to $Spec(k)$. A geometric point $\alpha$ in $S$ is a morphism from the spectrum $Spec(\bar{k})$ to $S$ where $\bar{k}$ is an algebraic closure/separable closure of $k$.} $Spa(\mathcal{C}) \rightarrow \mathcal{D}$
\cite{geo}. A profinite set is a compact totally disconnected topological space.%
\footnote{A profinite set $S$ is extremely disconnected if the closure of any open subset $U \rightarrow S$  is still open \cite{pro}.}$^{,}$%
\footnote{The connected components of a totally disconnected space are of the form $Spa(K, K+)$ for a perfectoid affinoid field $(K, K^+)$. For strictly totally disconnected spaces, one has in addition that K is algebraically closed.}  
A geometric point $Spa(\mathcal{C}) \rightarrow \mathcal{D}$ inside a diamond $\mathcal{D}$ is made ``visible'' by pulling it back through a quasi-pro-\'etale cover $X \rightarrow{\mathcal{D}}$, which gives rise to profinitely many copies of $Spa(\mathcal{C})$. 

While this property of perfectoid diamonds is taken to reflect the privacy of experience, talking about equivalence relations of classes of diamonds (and hence properties of ``selves'') needs more work. This is the scope of our three conjectures, focusing on $(\infty,1)$-categories of diamonds. In particular, what is needed is to prove dualizability and put a topology on such a category. This will eventually lead to the Efimov K-theory of large stable $(\infty,1)$-categories.

%

\subsection{Efimov K-theory of Diamonds $K^{Efimov} (\mathcal{Y}^{\diamond}_{S,E})$}

One of the authors has recently proposed to extend  K-theory  to the study of (categories of) perfectoid diamonds \cite{Dob21a}.
K-theory is defined on the category of small stable $\infty$-categories which are idempotent, complete and where morphisms are exact functors. A certain category of large compactly generated stable $\infty$-categories is equivalent to this small category. The condition for the colimit-preserving functor to preserve compact objects is that the morphisms are colimit preserving functors whose right adjoint also preserves colimits. 

In Efimov K-theory the idea is to weaken the condition of being compactly generated to being `dualizable' such that K-theory is still defined. A category $\mathcal{C}$ being dualizable implies that $\mathcal{C}$ fits into a localization sequence $\mathcal{C} \rightarrow \mathcal{S} \rightarrow \mathcal{X}$ with $\mathcal{S}$ and $\mathcal{X}$ compactly generated. Then the Efimov K-theory should be the fiber of the K-theory in the localization sequence.


\begin{figure}[htp]
   \centering
    \includegraphics[width=10cm]{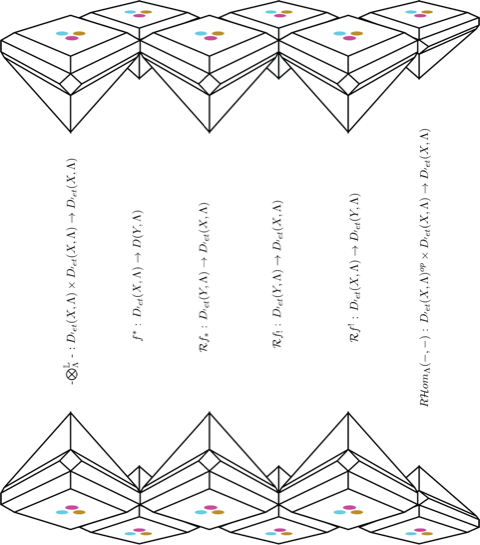}
   \caption{Efimov K-theory of diamonds as meta-model (picture taken from\cite{Dob21a}). While the perfectoid diamond is conjectured to model subjective experience, Efimov K-theory of diamonds describes the relations between classes of diamonds (as objects in an $(\infty,1)$-category).  Also shown are the definitions of the ``six operations'' (for details see text).
   }
  \label{fig:efimov-diamond}
\end{figure}
\noindent The technical details are the following:

\subsubsection*{Efimov Main Theorem}

Our terminology and small exposition follows Hoyois \cite{Hoy18}:

\begin{itemize}

\item $\infty$-categories are called categories. 

\item Let $\mathcal{P}r$ denote the category of presentable categories and colimit-preserving functors. 

\item Let $\mathcal{P}r^{dual} \subset \mathcal{P}r$ denote the subcategory of dualizable objects and right-adjointable morphisms (with respect to the symmetric monoidal and 2-categorical
structures of $\mathcal{P}r$). 

\item Let $\mathcal{P}r^{cg} \subset \mathcal{P}r$ be the subcategory of compactly generated categories and compact functors. Compact functors are functors whose right adjoints preserve filtered colimits.

\item Let $\mathcal{P}r^{\star}_{St}$ denote the corresponding full subcategories consisting of stable categories.

\item \textbf{Definition.} A functor $F : \mathcal{P}r^{dual}_{St} \rightarrow \mathcal{T}$ is called a localizing invariant if it preserves final objects and sends localization sequences to fiber sequences.

\item \textbf{Definition.} Let $C \in Pr$ be stable and dualizable. The continuous K-theory of $\mathcal{C}$ is the space
K$_{cont}(\mathcal{C}) = \Omega$ K$(Calk(C)^\omega)$.

\end{itemize}

\noindent \textbf{Lemma}. If $\mathcal{C}$ is compactly generated, then K$_{cont}(\mathcal{C}) = K(\mathcal{C}^w)$.\\

\noindent \textbf{Proof}. The localization sequence is Ind of the sequence $\mathcal{C}^w \hookrightarrow \mathcal{C} \rightarrow$ Calk$(\mathcal{C})^w$. Since K$(\mathcal{C}) = 0$, the result follows from the localization theorem in K-theory. \\

\noindent \textbf{Theorem [Efimov]}. Let $\mathcal{T}$ be a category. The functor 
\newline 
\indent $Fun(Pr^{dual}_{St} , \mathcal{T}) \longrightarrow Fun(Cat^{idem}_{St}, \mathcal{T})$, $F \mapsto F \circ Ind$,
\newline restricts to an isomorphism between the full subcategories of localizing invariants, with inverse $F \mapsto F_{cont}$.
In particular, if $\mathcal{C} \in Pr^{cg}_{St}$, then $F_{cont}(\mathcal{C}) = F(C^{\omega})$. \\

\noindent \textbf{Proof}. Let $\mathcal{A} \rightarrow \mathcal{B} \rightarrow \mathcal{C}$ be a localization sequence in $\mathcal{P}r^{dual}_{St}$. Then for large enough $\kappa$ we have an induced localization sequence
\newline 
\indent Calk$_{\kappa}(\mathcal{A}) \rightarrow$ Calk$_{\kappa}(\mathcal{B}) \rightarrow$ Calk$_{\kappa}(\mathcal{C})$.
\newline It follows that $F_{cont}$ is a localizing invariant. By Lemma we have $F_{cont} \circ$ Ind $\simeq F$, functorially in $F$. To $\mathcal{C} \in \mathcal{P}r^{dual}_{St}$ we can associate the filtered diagram of localization sequences $\mathcal{C} \rightarrow$ Ind$(\mathcal{C}^{\kappa}) \rightarrow$ Calk$_{\kappa}(\mathcal{C})$ for $\kappa \gg 0$. This gives a functorial isomorphism $F \simeq (F \circ$ Ind)$_{cont}$.\\

The main goal is to develop a Waldhausen S-construction to obtain the K-theory spectrum on the category of diamonds, and extend this to the $(\infty,1)$-category. 

In parallel, we will construct an $(\infty,1)$-site on the $(\infty,1)$-category of diamonds. To do this, we extend topological localization to the $(\infty,1)$-category of diamonds, recalling that equivalence classes of topological localizations are in bijection with Grothendieck topologies on $(\infty,1)$-categories $C$ \cite{top}. Topological localizations are appropriate to our construction because, under this type of localization of passing to the full reflective sub-$(\infty, 1)$-category, objects and morphisms have reflections in the category, mirroring the behavior of the perfectoid diamond. We give a small scale model of the localization procedure using Bousefield localization particular to triangulated categories.

Because a Grothendieck topology is a collection of morphisms designated as covers \cite{groth}, we find Scholze's six operations to be our basic building blocks in constructing the covers, in the following sense. Scholze's six operations live in the derived categories of sheaves and derived categories are the $\infty$-categorical localization of the category of chain complexes at the class of quasi-isomporhisms \cite{der}. The derived categories $\mathcal{D}(\mathcal{A})$ of abelian categories $\mathcal{A}$ are an important class of examples of triangulated categories. They are homotopy categories of stable $(\infty,1)$-categories of chain complexes in $\mathcal{A}$ \cite{triangulated}:

If we ignore higher morphisms, we can ``flatten" any $(\infty,1)$-category $\mathcal{C}$ into a 1-category ho(C) called its homotopy category. When $C$ is also stable, a triangulated structure captures the additional structure canonically existing on ho(C). This additional structure takes the form of an invertible suspension functor and a collection of sequences called distinguished triangles, which behave like shadows of homotopy (co)fibre sequences in stable $(\infty,1)$-categories.


\subsubsection*{$\mathcal{D}_{\acute{e}t}$ and the six operations}

The Grothendieck six operations formalism, also known as six functor formalism, is a formalization of aspects of 
the refinement of Poincar\'e duality from ordinary cohomology to abelian sheaf cohomology. 
This has been adapted by Scholze in the following way \cite{Sch17}:

\begin{itemize}
\item \textbf{Terminology}. Fix a prime $p$. Let $X$ be an analytic adic space on which $p$ is toplogically nilpotent. To $X$ was associate an 'etale site $X_{\acute{e}t}$. Let $\Lambda$ be a ring such that $n\Lambda = 0$ for some $n$ prime to $p$. There exists a left-completed derived category $\mathcal{D}_{\acute{e}t}(X, \Lambda)$
of \'etale sheaves of $\Lambda$-modules on $X_{\acute{e}t}$. Let $Perfd$ be the category of perfectoid spaces and $Perf$ be the subcategory of perfectoid spaces of characteristic $p$. Consider the $v$-toplogy on Perf \footnote{The $v$-topology, where a cover $\{f_i: X_i \rightarrow X\}$ consists of any maps $X_i \rightarrow X$ such that for any
quasicompact open subset $U \subset X$, there are finitely many indices $i$ and quasicompact open subsets
$U_i \subset X_i$ such that the $U_i$ jointly cover $U$ \cite{Sch17}}.
\item \textbf{Definition}.  A \textbf{diamond} is a pro-\'etale sheaf  $\mathcal{D}$ on $Perf$ which can be written as the quotient $X \slash R$ of a perfectoid space $X$ by a pro-\'etale equivalence relation $R \subset X \times X$. 
\\
Scholze restricts to a special class of diamonds which are better behaved. Denote the underlying topological space of a diamond
$Y$ as $|Y | $ = $|X| \slash |R|$.  
\item \textbf{Definition}. A diamond $Y$ is spatial if it is quasicompact and quasiseparated (qcqs), and $|Y|$ admits a basis for the topology given
by $|U|$, where $U \subset Y$ ranges over quasicompact open subdiamonds. More generally, $Y$ is locally
spatial if it admits an open cover by spatial diamonds.
\item \textbf{Definition} [from Definition 1.7 in \cite{Sch17}]. Let $X$ be a small $v$-stack, and consider the site $X_v$ of all perfectoid spaces over
$X$, with the $v$-topology. Define the full subcategory
$\mathcal{D}_{\acute{e}t}(X,\Lambda) \subset \mathcal{D}(X_v, \Lambda)$ as consisting of all $\mathcal{A} \in \mathcal{D}(X_v, \Lambda$ such that for all (equivalently, one surjective) map $f : Y \rightarrow X$
from a locally spatial diamond $Y$ , $f^*\mathcal{A}$ lies in $\hat{\mathcal{D}}(Y_{\acute{e}t},\Lambda)$.
\end{itemize}

\noindent $\mathcal{D}_{\acute{e}t}(X, \Lambda)$ contains the following six operations, see also Fig, \ref{fig:efimov-diamond}. 

\begin{enumerate}
\item Derived Tensor Product. -$\bigotimes^{\mathbb{L}}_{\Lambda}$ - $:D_{\acute{e}t}(X, \Lambda) \times D_{\acute{e}t}(X, \Lambda) \rightarrow D_{\acute{e}t}(X, \Lambda)$.

\item Internal Hom. $R\mathcal{H}om_{\Lambda}(-, -)$ : $D_{\acute{e}t}(X, \Lambda)^{op} \times D_{\acute{e}t}(X, \Lambda) \rightarrow D_{\acute{e}t}(X,\Lambda)$.

\item For any map $f : Y \rightarrow X$ of small $v$-stacks, a pullback functor $f^*$ : $D_{\acute{e}t}(X, \Lambda) \rightarrow D_{\acute{e}t}(Y, \Lambda)$.

\item For any map $f : Y \rightarrow X$ of small $v$-stacks, a pushforward functor $\mathcal{R}f_{*}$ : $D_{\acute{e}t} (Y, \Lambda) \rightarrow D_{\acute{e}t}(X, \Lambda)$.

\item For any map $f : Y \rightarrow X$ of small v-stacks that is compactifiable, representable in locally spatial diamonds, and with dim.trg $f < \infty$ functor $\mathcal{R}f_{!}$ : $D_{\acute{e}t} (Y, \Lambda) \rightarrow D_{\acute{e}t}(X, \Lambda)$.

\item For any map $f : Y \rightarrow X$ of small $v$-stacks that is compactifiable, representable in locally spatial
diamonds, and with dim.trg $f < \infty$, a functor $\mathcal{R}f^{!}$ : $D_{\acute{e}t} (X, \Lambda) \rightarrow D_{\acute{e}t}(Y, \Lambda)$.

\end{enumerate}

%

\subsection{Organization of the paper}

The notion of ``1-awareness'' refers to a network of representations, for example, when concentrating on an exam or a zoom call, while being implicitly directed to the environment and to future events.  
If a category theoretical approach is correct, it is natural to also consider higher representational levels, based on 1-awareness. We thus refer to ``2-awareness'' as  simultaneously sustaining any combination of 1-awarenesses, such as in the following example: Imagine you are having lunch in your home on a Saturday at noon in the year 2021 while simultaneously you are a Cambridge Apostle having a discussion in the Moral Sciences Club on a Saturday evening in 1888. Continuing, 3-awareness would be the sustaining of any combination of 2-awarenesses. For example, simultaneously being aware of simultaneously being aware of having lunch in your home on a Saturday at noon localized in the year 2021 while being a Cambridge Apostle having a discussion in the Moral Sciences Club on a Saturday evening in 1884, and of giving a thesis defense on perfectoid spaces on a Friday at 10am in 2024, while being at your desk and finishing the first chapter of said thesis later this evening... 
And so on, up to  ``$n$-awareness''.

In section \ref{sec:model}, we model this idea based on higher category theory (the ``category of $n$-awareness'', starting with $n=1$). The objects of ``1-awareness'' are representations; morphisms between such objects are maps between such representations. 

We then distinguish ``higher'' categories, based on the category of 1-awareness. The set $Hom(X,Y)$ is the set of all morphisms between objects $X$ and $Y$ in a category. We would further want $Hom(X,Y)$ to be more than a set. We would like $Hom(X,Y)$ to be a topological space of morphisms from $X$ to $Y$.  
%
One can appeal to homotopy theory to precisely study transformation and equivalences between such topologically structured spaces.\footnote{In order to describe this more precisely, we will, in Conjecture 1, conjecture that the complexities of covering maps can be represented in terms of the pro-\'etale topology on the category of perfectoid spaces.} 
We conclude this section by generalizing our setup to $(\infty,1)$-topoi to look at objects with higher homotopy. This defines our basic meta-model for $n$-awareness (with $n\rightarrow\infty$ for $(\infty,1)$-topoi).

In Section \ref{sec:conjectures} we present three conjectures to ``enrich'' our framework. First, we propose that subjective experience can be modeled in terms of a geometric structure known as ``perfectoid diamond'' \cite{SW20}. Specifically, we first propose that the  $(\infty,1)$-category of perfectoid diamonds is an $(\infty,1)$-topos. 
To study the topology on such $(\infty,1)$-categories,
we, secondly, propose that topological localization extends to the $(\infty,1)$-category of diamonds. 
We provide a small-scale model, using triangulated categories. Triangulated categories are homotopy categories of stable $(\infty,1)$-categories. 
Thirdly, Efimov K-theory allows one to study localization sequences of sheaves on such categories. 
Together, these three conjectures define our meta-model, satisfying the constraints embodied within our framework and inheriting several properties such as the descent condition of $(\infty,1)$-topoi.

In section \ref{sec:consequences}, we discuss implications of this meta-model. Those pertain to language, time, and selfhood. 
%
%
We first sketch the grammar (``$n$-declension'') of a new language, which follows from our mathematical model of $n$-awareness and its concomitant extension to $n$-time. 
Such a linguistic system enables us to actually express experience, revealing the way  we concretely engage with the world -- and constantly make the error of suspecting substantial (metaphysical) entities where there are none, analogous to the conception of ``language games'' \cite{Wit53}, updated to a ``language scheme'' inspired by the mathematical model proposed.

The fiction of substantial temporal parts is rejected. How to regard subjective experience against (the ontology of) time is one of the most puzzling open questions in philosophy. And any mathematical treatment of subjective experience should be able to shed some light on its relation to time.  \textit{Time} in our framework is modeled as hierarchy of  ``$k$-times'' ($k\le n$): 1-awareness implies a notion of 1-time, 2-awareness implies a notion of 2-time, ... , and $n$-awareness implies a notion of $n$-time. On our reading, temporality is just a consequence of the geometry of representation. It is at most an ordering property. This is not an altogether new idea \cite{Hug18}, but our mathematical model bears some natural consequences that have not, to our knowledge, been developed so far in the literature: our structural definition of $n$-time puts a novel spin on the problem of temporal coincidences and simultaneity. A 2-awareness contains its 1-awarenesses, which are thus ``simultaneously present'' in 2-awareness. Moreover, we shall explain how this simultaneity could be understood within our meta-model.

We then turn to the concept of ``\textit{selfhood}". One desideratum that any theory of selfhood should meet is that any self's experience is self-reflexive, i.e. all moments (parts) of such an experience are referring to the total experiential context of the self. This distinguishes the experience of a self from experience ``as such''.  We account for this self-reflexivity by the idea of having invertible morphisms between objects in higher categories, implying a ``mirroring'' of representations in all other representations -- an ``internalized monadology'' in the spirit of Leibniz \cite{Lei14}. This attributes self-reflexivity on the level of the next higher category. Since $(k-1)$-categories are contained as objects in a $k$-category, the property of self-reflexivity is ``lifted'' to  the next higher level of awareness.
However, having strictly invertible morphisms would imply structural equivalence. 
This would destroy any notion of a ``true individual'' on the level of selves. The philosophical problem underlying this dilemma is the following: If the universe evolves toward a system of equivalent representations, how to prevent it eventually becoming ``totalitarian'', neglecting individuality? However, in the framework of homotopy theory, there is a sense in which such equivalences only hold 
%
``up to'' \cite{upt} homotopy equivalence.%
\footnote{We later, in Conjecture 2, conjecture to work with abelian groups \cite{abe} and to form a complex \cite{com} of these groups, which is informally a sequence of compositions. Then we  ``chain'' the complexes together using a chain map. These complexes are the objects in the ``derived category'' \cite{Art63}. Thus we have upgraded representations from being objects in a general category to being objects in the derived category.}
%
%
So, rather than thinking about ``selves'' as enduring entities wholly present at every moment in time (i.e. like endurantists would), but also unlike the idea that we have (more or less well-defined) ``temporal parts'' (i.e. like perdurantists would), we replace the idea of personal identity over time with an updated form of the ``bundle theory of self''\footnote{The ``bundle'' refers to a collection of (weakly equivalent) awareness-objects which are related  via morphisms.} \cite{Hum40}, based on a notion of ``self-reflexivity'' that preserves the individuality of experience. 
\\
Some additional mathematical definitions can be found in an appendix.
\section{The vanilla framework}
\label{sec:model}
\subsection{Basics of higher category theory}

A mathematical model should express things in a minimal and conceptually coherent way. (Higher) category theory, while seemingly abstract, offers a flexible framework in which to define ``relations'', arguably the basic structural property of subjective experience \cite{Ehr15,Tsu16,Nor19,Kle20,Sig20}. 
Our meta-model (``$n$-awareness'') should not be regarded as model of a substantial (e.g. spatio-temporal) entity, but as relational object -- as if the structure alone dictated its properties. In particular, the following bears some similarities to the work of Andr\'ee Ehresmann and coworkers \cite{Ehr15,Ehr07}, however with the difference that the primitive objects are not physiological entities (e.g. neurons), but representations.

In mathematics, category theory formalizes the notion of ``structure'' in a very general way. A category is defined by its \textit{objects}, $A,B$, and  \textit{morphisms} (structure-preserving maps) between those objects, $A \rightarrow B$, satisfying certain requirements (composing associatively and the existence of an identity).\footnote{Two examples are the category of sets where this amounts to the fact that morphisms are functions, or the symmetric monoidal category where this condition means to satisfy hexagonal identity.} 
We propose to change the ``category number'' of mathematical models to study subjective experience. Awareness is sometimes identified with a ``0-category'', e.g. computable functions between domains.
One prominent approach along these lines is machine-state functionalism \cite{Put67} which identifies ``mental states'' with functional relations involving sensory inputs and behavioral outputs, representable as string of numbers. Importantly, functionalism also assumes that mental states are individuated with respect to the relations they have to other mental states.

Upgrading this idea further, we propose ``$n$-categories'' to study subjective experience:  A ``$1$-category'' contains objects and ``1-morphisms'' which are morphisms between these objects. 
For disambiguation, the objects of the categories could be referred to as ``moments of awareness'', defined as dependent parts within an experiential whole, following an early idea of Edmund Husserl \cite{Hus75}
-- in the case of 1-awareness we call those objects ``representations'' (Table \ref{tab:cat}).

\begin{table}[hbt]
\caption{Schematic for understanding subjective experience in terms of a (nested) hierarchies of categories (``$n$-awareness'').}
\vspace*{0.25cm}
\begin{tabular}{c| l}
\toprule
& $n$-awareness: a category of categories of... representations\\
$\uparrow$& $\vdots$\\
subjective experience & 3-awareness: a category of categories of categories of representations \\
$\downarrow$ & 2-awareness: a category of categories of representations \\
& 1-awareness: a category of representations \\
\bottomrule
\end{tabular}
\label{tab:cat}
\end{table}

We do not wish to enter the debate surrounding the question \textit{what}\footnote{Some popular proposals: functional states \cite{Blo96}, patterns of brain states \cite{Wu18}, ``intentional'' objects of thought \cite{Cra13}, ``monads'' \cite{Lei14}, or fundamental constituents of the universe \cite{Hof08,Hof14,Gof17}.}
 these representations are, but we are primarily interested in \textit{how they relate}.  Our hope is that we can model certain interesting properties of subjective experience by establishing a framework for modeling the relations between representations  (and how this is complexified). 
We hence do not start with a metaphysical ``working definition'' of what subjective experience is, but give a model of what it \textit{does}.

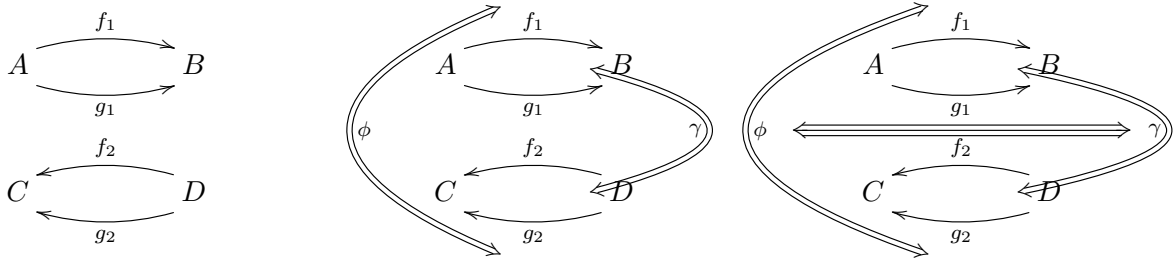
\begin{figure}
\begin{equation}
\nonumber
\xymatrixrowsep{0.5cm}
\xymatrixcolsep{.8cm}
\xymatrix{
&&&&&&\ar@2{<->}@<-2ex>@/_5pc/[dddd]^\phi &&&&&\ar@2{<->}@<-2ex>@/_6pc/[dddd]^\phi&\\
A \ar@/^/@<1ex>[rr]^{f_1} \ar@/_/@<-1ex>[rr]_{g_1} && B &&  &
A \ar@/^/@<1ex>[rr]^{f_1} \ar@/_/@<-1ex>[rr]_{g_1} & \ar@2{<->}@<4ex>@/^4pc/[dd]_\gamma & B && &
A \ar@/^/@<1ex>[rr]^{f_1} \ar@/_/@<-1ex>[rr]_{g_1} & \ar@2{<->}@<4ex>@/^5pc/[dd]_\gamma& B && &
\\
&&&&  &  &&& &\ar@3{<->}[rrrr] &&&&\\
C && \ar@/^/@<1ex>[ll]^{g_2} \ar@/_/@<-1ex>[ll]_{f_2} D &&  &
C & & \ar@/^/@<1ex>[ll]^{g_2} \ar@/_/@<-1ex>[ll]_{f_2} D && &
C & & \ar@/^/@<1ex>[ll]^{g_2} \ar@/_/@<-1ex>[ll]_{f_2} D\\
&&&&&&&&&&&&
} 
\end{equation}
\caption{A graphical overview for $n$-awareness, for $n=1,2,3$. Left: 1-category consisting of objects $A,B,C,D$ (representations) and possible 1-morphisms, $f_i, g_i$ between $A$ and $B$ (respectively between $D$ and $C$). Middle: 2-category which also includes 2-morphisms,  $\phi: { f_1 \atop g_1} \rightarrow {f_2 \atop g_2}$  and $\gamma:  { f_1 \atop g_1} \rightarrow {g_2 \atop f_2}$. Right: 3-category which also includes 3-morphisms between 2-morphisms.}
\label{fig:n-cat}
\end{figure}

1-awareness consists of various representations and the relations between them. This means that, for example, me being aware of the zoom call right now is contained (as representation) within a category populated by myriads of other (also potential) moments of awareness -- sensing the joy in your voice, or being afraid of your dismissive reactions to what I want to say -- but also my (implict) perceptions of environmental going-ons, memories, and other (possibly explicit) background experiences; and similar for the examples of attending a meeting of the Cambridge apostles in 1888 (this experience is related to my potential experience of listening to Wittgenstein a few years later), or for the experience of a thesis defense in 2024 (this experience is related to my actual experiences of writing the thesis now), etc. See Fig. \ref{fig:n-cat} (left).

One could ``increase the category number'' to study possible relations between relations. A ``$2$-category'' contains objects, $1$-morphisms, and $2$-morphisms.  In Fig. \ref{fig:n-cat} (middle) we not only have $1$-morphisms between objects but also $2$-morphisms between the $1$-morphisms. The collection of all morphisms from $A$ to $B$ forms a set called  the ``homset'', $Hom(A,B)$. In a 2-category, each homset itself carries the structure of a category -- a collection of objects and morphism satisfying certain requirements -- and thus morphisms between such homsets can be regarded as morphisms between 1-categories. 
This higher dimensional structure allows the 2-category to sustain two moments of awareness ``at once''.  See Fig. \ref{fig:n-cat}, (middle).
Continuing further, a ``3-category'' contains objects, $1$-morphisms, $2$-morphisms, and $3$-morphisms between the $2$-morphisms, where $3$-morphisms can be seen as relations between 2-categories (or the respective moments).  (Fig. \ref{fig:n-cat}, right). 

It follows that an $n$-category contains objects, $1$-morphisms, $2$-morphisms, …, up to $n$-morphisms between the $(n-1)$-morphisms.  Just as $2$-awareness ``contains'' $1$-awareness, this structure is nested. 
%

\subsection{Homotopy theory}

Instead of having just a \textit{set} of morphisms from one object to another (i.e. the homset), we want our categories to have the structure of a \textit{topological space} of morphisms from one object to another\footnote{Others have previously hypothesized, for slightly differrent reasons, that a mathematical model of consciousness should have at least the structure of a topological space \cite{Pre19,Kle20}.}
, so that we can have a \textit{space} of maps between the objects we wish to study. This also allows us to speak about transformations and equivalences between such spaces, and thus between moments of awareness.  
We appeal to
\textit{homotopy theory} to make these ideas more precise. %
This will allow us to understand $n$-awareness as (nested) ``simultaneous presence'' of $m$ moments of 1-awareness, thereby avoiding certain problems that are related to time. Homotopy theory, which studies deformation equivalences called homotopies, is defined as follows: Two continuous functions from one topological space to another are called homotopic if one can be continuously deformed into the other. \textit{Homotopy groups} extend this notion to equivalences between topological spaces, to therein classify topological spaces.\footnote{We say that two topological spaces, $X$ and $Y$, are of the same homotopy type or are homotopy equivalent if we can find continuous maps $f: X \rightarrow Y$ and $g: Y \rightarrow X$ such that $g \circ f$ is homotopic to the identity map $idX$ and $f \circ g$ is homotopic to $idY$.} There is a remarkable freedom in reducing strong equivalences, such as the claim that object $A$ equals object $B$, to deformation equivalences, such as the claim that object $A$ is ``deformation equivalent'' to object $B$ because one can create the homotopy map which makes these objects homotopy equivalent. This allows us to speak of two experiences as being equivalent ``up to'' homotopy.
\subsection{Infinity Topoi}

A topos (Greek for ``place'') is a category which behaves like the category of sets but also contains a notion of localization.%
\footnote{Topoi are modeled after Grothendieck’s notion of a sheaf on a site \cite{Luc04}. A prototypical example of a topos is the category of sets, since it is the category of sheaves of sets on the one point space. Informally, topoi are ``nice'' categories for doing geometry that act like models of intuitionist type theory. They are abstract contexts ``in which one can do mathematics independently of their interpretation as categories of spaces'' \cite{Art63,top}.}
%
Morphisms between objects of $n$-awareness represent the way how different moments of awareness are related to each other. In case morphisms are \textit{invertible}, a mapping from $Hom(A,B)$ to $Hom(C,D)$ implies the existence of a mapping from $Hom(C,D)$ to $Hom(A,B)$.
An $(\infty,k)$-category is an infinity category ($n\rightarrow\infty$) in which all morphisms higher than $k$ are invertible. For our model, we chose to work in $(\infty, 1)$-categories, where all $k$-morphisms (for $k>1$) are invertible; this means that the respective moments of awareness could be transformed into an equivalent structure. This also means that any moment of $k$-awareness ``reflects'' all the other objects it is related to, and looking at the whole category (or at the respective object in the next higher category), it can be said that $n$-awareness is ``self-reflexive''.

An $(\infty,1)$-topos is a generalization of topos theory to higher category theory.\footnote{In this context, higher category theory investigates the generalizations of $\infty$-groupoids to directed spaces \cite{hig}.} The objects in $(\infty,1)$-topoi are generalized spaces with higher homotopies that carry more structure. A prototypical example would be $(\infty, 1)$-topoi of $(\infty,1)$-sheaves \cite{inf}.
More precisely,  an $(\infty, 1)$-topos is an $(\infty,1)$-category $C$ which satisfies three conditions: $C$ is presentable (with kappa filtered colimits), locally cartesian closed, and satisfies a descent condition (where an object in $C$ is sent to the slice category $C/u$). We propose these conditions to reflect the property of subjective experience to be ``self-contained''.
\footnote{In general, topoi are models of internal logic, which means that almost any \cite{top} 
logical statement could be internalized. $(\infty,1)$-topoi are types with internalized descent datum \cite{dec}.}
All conceivable relations are built inside the framework. 

\begin{figure}[htb]
\begin{equation}
\nonumber
\xymatrixrowsep{0.4cm}
\xymatrixcolsep{.6cm}
\xymatrix{
& K^{Efimov}(\mathcal{Y}^{\diamond}_{S,E}) &\\\\\\\\
\txt{$\infty$-categories \\ + invertible $k$-morphisms  $(k>1)$} \ar[dr]    &  \includegraphics[width=4cm]{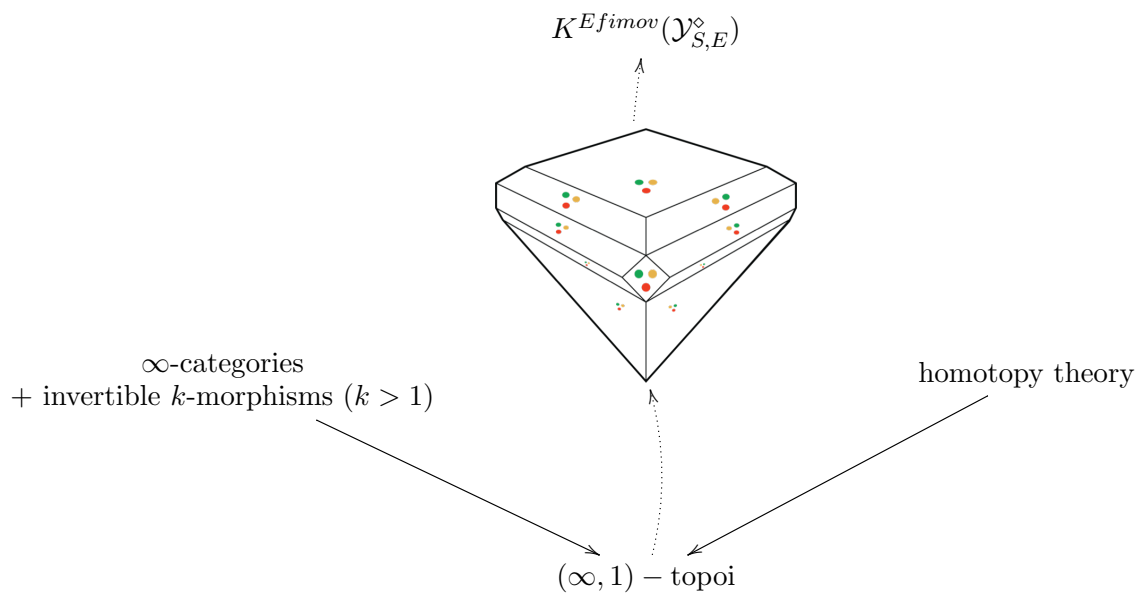} \ar@{.>}@/^/[uuuu] && \textrm{homotopy theory}\ar[dll] \\
& (\infty,1)-\textrm{topoi}\ar@{.>}@/_/[u]
}
\end{equation}
\caption{\textit{Gestell} of a meta-model. 
The theory of $(\infty,1)$-topoi combines higher category theory with the notions of homotopy equivalence between topological spaces and invertible morphisms. The  $(\infty,1)$-category of perfectoid diamonds is conjectured to be an $(\infty,1)$-topos (Conjecture 1), and its structural equivalences to subjective experience can be made precise within Efimov K-theory (Conjecture 3).}
\label{fig:mode}
\end{figure}

The $(\infty,1)$-topos is, informally, a place where one can do homotopy theory, and it represents our meta-model of subjective experience.  
The framework is visualized in Fig. \ref{fig:mode}, together with a graphical overview of the conjectured meta-model that is described in section \ref{sec:conjectures}. Also shown is a conjectured Efimov representation of the form $\mathcal{Y}^{\diamond}_{S,E}$ = $S \times (Spa \mathcal{O}_E)^{\diamond}$, the diamond Fargues-Fontaine curve of the geometrization of the local Langlands correspondence \cite{Far16}.

\section{Arriving at the meta-model: Three conjectures}
\label{sec:conjectures}

The vanilla framework of $(\infty,1)$-topoi needs to be enriched to capture more properties typically associated with subjective experience. We now wish to propose three conjectures for the mathematical study of subjective experience.
%
\\\\
\textbf{Conjecture 1}. The  $(\infty,1)$-category of perfectoid diamonds is an $(\infty,1)$-topos.
\\\\
\textbf{Remark 1}. Recall, the category of sheaves on a (small) site is a Grothendieck topos.
\\\\
\textbf{Remark 2}. Lurie discusses the structure needed for this construction. Recall the following:
\newline \\
\textbf{Definition} [see also Definition 6.2.2.6 in \cite{lur}]: 
An $(\infty,1)$-category of $(\infty,1)$-sheaves is a reflective sub-$(\infty,1)$-category $Sh(C) \overset{\overset{L}{\leftarrow}} {\hookrightarrow} PSh(C)$
of an $(\infty,1)$-category of $(\infty,1)$-presheaves such that the following equivalent conditions hold:

\begin{itemize}
\item 1. $L$ is a topological localization.
\item 2. There is the structure of an $(\infty,1)$-site on $C$ such that the objects of $Sh(C)$ are precisely those $(\infty,1)$-presheaves $A$ that are local objects with respect to the covering monomorphisms $p:U \rightarrow j(c)$ in $PSh(C)$ in that

\hspace*{0.25cm} $A(c) \simeq PSh(j(c),A)$ $\xrightarrow{\text{PSh(p,A)}}$ $PSh(U,A)$
is an $(\infty,1)$-equivalence in $\infty$Grpd.
\item 3. The $(\infty,1)$-equivalence is the descent condition and the presheaves satisfying it are the $(\infty,1)$-sheaves.
\end{itemize}

%
\noindent \textbf{Remark 3}. This conjecture requires enlarging the class of higher pro-\'etale morphisms of diamonds, restricting to locally spatial diamonds, and is motivated in \cite{SW20}:
\begin{quote}
Theorem 9.1.3 (\cite[Propositions 9.3, 9.6, 9.7]{Sch17}). Descent along a pro-\'etale cover $X' \rightarrow X$ of a perfectoid space $f : Y' \rightarrow X'$
is effective in the following cases.
\begin{itemize}
\item 1. If $X, X'$ and $Y, Y'$ are affinoid and $X$ is totally disconnected.
\item 2. If $f$ is separated and pro-\'etale and $X$ is strictly totally disconnected.
\item 3. If $f$ is separated and \'etale.
\item 4. If $f$ is finite \'etale.
\end{itemize}
Moreover, the descended morphism has the same properties.
\end{quote}

\noindent The diamond properties of coherence (quasicompact, and quasiseparated), the pro-\'etale topology, and its presentability correlate with assumed properties of subjective experience. In a perfectoid diamond, privacy is mirrored in its interior ``geometric points'' being made visible as profintely many covers. Geometric points have different appearances depending on the ``angle'' they are viewed from, and complete information about these points cannot be recovered.
Informally, the structure of the diamond is defined topogically (by its pro-\'etale topology). Self-containedness is a property inherited from the $(\infty,1)$-topos of which the $(\infty,1)$ category of diamonds is a conjectured example.\\


\noindent \textbf{Definition}. Let $Perfd$ be the category of perfectoid spaces and $Perf$ be the subcategory of perfectoid spaces of characteristic $p$. 
A \textbf{diamond} is a pro-\'etale sheaf  $\mathcal{D}$ on $Perf$ which can be written as the quotient $X \slash R$ of a perfectoid space $X$ by a pro-’etale equivalence relation $R \subset X \times X$. 
\\
In essence, a diamond is a sheaf $Y$ on the category Perf which is constructed as the quotient of a perfectoid space by a pro-\'etale equivalence relation. More formally, this diamond quotient lives in a category of sheaves on the site of perfectoid spaces with pro-\'etale covers.\\

\noindent One \textbf{example} of a diamond is given by the sheaf $SpdQ_p$ that attaches to any perfectoid space $S$ of characteristic $p$ the set of all untilts $S \#$ over $Q_p$, where $\textit{untilting}$, informally, is an operation that translates from characteristic $p$ back to characteristic $0$ $\cite{Sch12}$. The formal definition is:
\\\\
\textbf{Definition}. $SpdQ_p = Spa(Q_p^{cycl})/ \underline{Z_p^\times} $. That is, $SpdQ_p$ is the coequalizer of $\underline{Z_p^\times} \times Spa(Q_p^{cycl})^\flat \rightrightarrows Spa(Q_p^{cycl})^\flat$,  where one map is the projection and the other is the action; cf. \cite{SW20}.\\

\noindent A second \textbf{example} is the moduli space of shtukas for $(\mathcal{G}, b, \{ \mu_1, . . . , \mu_m) \}$ fibered over the m-fold product $SpaQ_p \times SpaQ_p ...\times_m SpaQ_p$ \cite{SW20}.
\\

\noindent Perfectoid spaces provide a framework for translating number-theoretic class field problems from characteristic $p$ to characteristic $0$. The definition of perfectoid space is as follows.
\\\\
\textbf{Definition}. A perfectoid space is an adic space covered by affinoid adic spaces of the form $Spa(R, R^+)$ where $R$ is a perfectoid ring. 
\\\\
\textbf{Examples} of perfectoid spaces, as detailed in \cite{SW20}, are the following: 

\begin{itemize}
\item The Lubin-Tate tower at infinite level: $\mathcal{M}_{LT,\infty} = \tilde{U}_x \times^{GL_2(\mathbb Q_P)_1}GL_2(\mathbb Q_P)\cong \underset{Z}{\sqcup} \tilde{U}_x$ \cite{SW20}. 

\item Perfectoid Shimura varieties: $S_{K^p} \sim \lim\limits_{\overleftarrow {K^p}}(S_{K^pK_p} \bigotimes_E E_p)^{ad}$ \cite{Sch12}.

\item Any completion of an arithmetically profinite extension (APF), in the sense of Fontaine and Wintenberger {\cite{FW79}}, is perfectoid. A nice source of APF extensions is $p$-divisible formal group laws. 

\item If $K$ is a perfectoid field and $K+ \subset K$ is a ring of integral elements, then $Spa(K, K+)$ is a perfectoid space.

\item Zariski closed subsets of an affinoid perfectoid space support a unique perfectoid structure. 

\item It is mentioned that the nonarchimedean field $Q_p$ is not perfectoid as $Q_p$ does not have a topologically nilpotent element $\xi \in Z_p$ whose $p$-th power divides $p$.
\item \textbf{Definition.} A perfectoid space $X$ is totally disconnected if it is quasicompact and quasiseparated (qcqs) and every open cover splits. A perfectoid space $X$ is strictly totally disconnected if it is qcqs and every \'etale cover splits. Any totally disconnected perfectoid space is also affinoid.
\end{itemize}

%

Why do we use adic spaces? Informally, adic spaces are versions of schemes associated to certain topological rings and it is these structures that produce the nonarchimedean fields we need to represent the local/global nature of consciousness. Formally, an adic space is a triple $(X, \mathcal{O}_X,\mathcal{O}_X^+)$ where $X$ is a topological space, $\mathcal{O}_X$ is a sheaf of complete topological rings, and $\mathcal{O}_X^+$is a subsheaf of $\mathcal{O}_X$ that is locally of the form $(Spa(A, A^+), \mathcal{O}_A, \mathcal{O}_A^+ )$. 

We assume that conscious agents operate via topological coverings, i.e. via establishing the right clustering of (data) points, rather than via calculating a particular function based on such a clustering. But in order to describe this process, we need a way to capture the complexities of covering maps
that preserves local homeomorphism and diffeomorphism properties. \'Etale morphisms encapsulate the idea of maps which are local homeomorphisms/local diffeomorphism. 
The pro-\'etale site is the site of perfectoid spaces with pro-\'etale covers.
\footnote{The pro-\'etale site tempers the ``finiteness conditions" on the fibers of \'etale morphisms to pro-finiteness conditions, which are pro-\'etale morphisms; cf. Remark 10.3.2 in \cite{SW20}.}


Moreover, the language of sheaves provides a geometrization that can model the  global storage of data and recollection of local data. Working over structures like $Spec R$, the spectrum of a ring $R$, allows us to add additional structure to points, encoding geometric parameter spaces or  moduli spaces. A more sophisticated mathematical structure allows us to model more complex behavior. Subjective experience seems to proceed from a local to global setting and back, mimicking the tilting and untilting operation of perfectoid spaces from characteristic $0$ to characteristic $p$.
\\\\
\textbf{Conjecture 2}.  Topological localization, in the sense of Grothendieck-Rezk-Lurie $(\infty,1)$-topoi, extends to the $(\infty, 1)$-category of diamonds.
\\\\
\textbf{Remark 4}. To construct a topology on the $(\infty,1)$-category of diamonds, we must first construct the $(\infty,1)$-site on the $(\infty,1)$-category of diamonds. Recall, the definition of an $(\infty,1)$-site.
\\
\newline \textbf{Definition} \cite{site}. The structure of an $(\infty,1)$-site on an $(\infty,1)$-category $\mathcal{C}$ is precisely the data encoding an $(\infty,1)$-category of $(\infty,1)$-sheaves $Sh(\mathcal{C}) \hookrightarrow PSh(\mathcal{C})$
inside the $(\infty,1)$-category of $(\infty,1)$-presheaves on $\mathcal{C}$.
\\
\newline Topological localizations are appropos to our diamond construction as objects and morphisms have reflections in the category, just as geometric points have reflections in the profinitely many copies of $Spa(\mathcal{C})$. Recall the definition of a reflective subcategory. 
\\
\newline \textbf{Definition} \cite{reflective}. A reflective subcategory is a full subcategory
$\mathcal{C} \hookrightarrow \mathcal{D}$
such that objects $d$ and morphisms $f:d \rightarrow d'$ in $\mathcal{D}$ have ``reflections'' $Td$ and $Tf:Td \rightarrow Td'$ in $\mathcal{C}$. Every object in $\mathcal{D}$ looks at its own reflection via a morphism $d \rightarrow Td$ and the reflection of an object $c \in \mathcal{C}$ is equipped with an isomorphism $Tc \simeq c$. The inclusion creates all limits of $\mathcal{D}$ and $\mathcal{C}$ has all colimits which $\mathcal{D}$ admits.
\\\\
\textbf{Remark 5}. Extending a perfectoid version of localization is essential to the construction of the $(\infty,1)$-topoi. See \cite{SW20}, Remark 10.3.2. This is not straightforward as the topology of the diamond $\mathcal{D}$ can be highly pathological and not even $T_0$\footnote{ An example from \cite{SW20} details the quotient of the constant perfectoid space $Z_p$ over a perfectoid field
by the equivalence relation “congruence modulo $Z$” which yields a diamond with
underlying topological space $Z_p / Z$.}. However, we will restrict to a special class of well-behaved diamonds, qcqs, given in \cite{SW20}: 
\begin{quote}
Proposition 10.3.4. Let $\mathcal{D}$ be a quasicompact quasiseparated diamond. Then $\mathcal{D}$ is $T_0$.
\end{quote}

%
\noindent We previously appealed to the ``up to'' notion that comes naturally with homotopy theory and to $k$-morphisms being invertible to mirror the self-reflexivity property. 
To make the idea of the associated equivalence (``equivalent up to homtopy'') more precise, we invoke the idea of a \textit{weak equivalence}.  We use derived category theory \cite{Art63} to form a complex \cite{com} of abelian groups, which is informally a sequence of compositions. Then we  ``chain'' the complexes together using a chain map. Thus we have upgraded representations from being objects in a general category to being objects in the derived category. Weak equivalences can also be referred to as ``quasi-isomorphisms'' in this context.

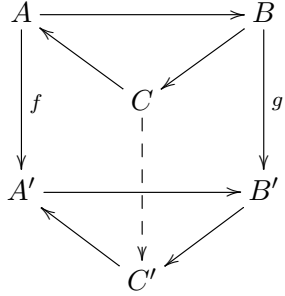
\begin{figure}
\begin{equation}
\nonumber
\xymatrixrowsep{.7cm}
\xymatrixcolsep{1.cm}
\xymatrix{
A \ar[rr]  \ar[dd]^{f} && B \ar[dl] \ar[dd]^{g} \\
& C \ar[ul] \ar@{-->}[dd] & \\
A' \ar[rr] && B' \ar[dl]\\
& C' \ar[ul] &
}
\end{equation}
\caption{A triangulated category composed of two exact triangles $A, B, C$ and $A’ , B’, C’$ and three commuting squares $ABB’A’$, $AC'C'A'$, and $BB'C'C$, with a non-unique fill-in \cite{Cal05}.}
\label{fig:triang}
\end{figure}

Let $\mathcal{A}$ be a Grothendieck abelian category (e.g., the category of abelian groups). We define $K(\mathcal{A})$ to be the homotopy category of $\mathcal{A}$ whose objects are complexes of objects of $\mathcal{A}$ and whose homomorphisms are chain maps modulo homotopy equivalence. 
The weak equivalences are quasi-isomorphisms defined as follows: A chain map $f: X \rightarrow Y$ is a quasi-isomorphism if the induced homomorphism on homology is an isomorphism for all integers $n$. 
$K(\mathcal{A})$ is endowed with the structure of a \textit{triangulated category}. A triangulated category has a translation functor and a class of exact triangles which generalize fiber sequences and short exact sequences. Localization by quasi-isomorphisms preserves this triangulated structure. Bousfield localization \cite{bou}\footnote{
Through this process, the derived category $D(\mathcal{A})$ of the initial abelian category $\mathcal{A}$ is obtained by ``pretending” that quasi-isomorphisms in $K(\mathcal{A})$ are isomorphisms. 
Specifically, the localization is constructed as follows: morphisms in $D(\mathcal{A})$ between $A\cdot$ and $B\cdot$ will be 'roofs' \cite{Cal05}, with $f, g$ morphisms in $K(\mathcal{A})$ and $f$ a quasi-isomorphism. This roof represents $g \circ f^{-1}$.}, 
particular to triangluated categories, allows us to make more morphisms count as weak equivalences and this is formally how we get from 1-awareness to $n$-awareness. One of the axioms of a triangulated category states that given the diagram in Fig. \ref{fig:triang},
where $A, B, C$ and $A’ , B’, C’$ form exact triangles, and the morphisms $f$ and $g$ are given such that the square $ABB’A’$ commutes, then there exists a map $C \rightarrow C’$ such that all the squares commute. This triangulated category is a categorization of a set-theoretic \textit{ordinal}. A set is an ordinal number if it is transitive and well-ordered by membership, where a set $T$ is transitive if every element of $T$ is a subset of $T$.  Ordinals locate within a set as opposed to  cardinality which references merely size.


$K(\mathcal{A})$ represents all possible moments of awareness (related by composition). There are an infinite number of 1-awarenesses in this category. Differentials define a 1-awareness relating with another 1-awareness. The composition of any two $A$’s is the zero map. This is the group law, the return to identity. Let $B$ be the abelian group of 2-awarenesses.  The differentials define relating with another 2-awareness. The diagram commuting means the 2-awareness is related with the 1-awareness. 
In essence, this diagram represents a higher-dimensional notion of commutativity through the map $C \rightarrow C'$ and the three squares commuting. By higher-dimensional, we mean the three commuting squares and the two exact triangles together form a cone construction. There is a relation here between the higher commutativity and extending the pro-\'etale morphisms of diamonds.\\\\
\noindent \textbf{Conjecture 3}. The meta-model takes the form of Efimov K-theory of the large stable $(\infty,1)$-category of perfectoid diamonds. 
\\\\
\noindent This will require proving the dualizability of the category, see \textbf{Theorem [Efimov]} \cite{Dob21a}.
Specifically, we seek a diamond reformulation of the following Efimov theorem \cite{Dob21a}.\\

\noindent \textbf{Theorem [Efimov*]}. Let $X$ be a locally compact Hausdorff topological space, $\mathcal{C}$ a stable dualizable presentable category, and $R$ a sheaf of E${_1}$-ring spectra on $X$. 
Suppose that $Shv(X)$ is hypercomplete (i.e., $X$ is a topological manifold). Let $\mathcal{T}$ be a stable compactly generated category and $F $: Cat$^{idem}_{St} \rightarrow \mathcal{T}$ a localizing invariant that preserves
 filtered colimits. Then:\\

$F_{cont}(Mod_R(Shv(X, \mathcal{C}) \simeq \Gamma_c(X, F_{cont}(Mod_R(\mathcal{C})))$.
Specifically, $F_{cont}(Shv(\mathbb{R}^n, \mathcal{C})) \simeq \Omega^nF_{cont}(\mathcal{C})$.
\\\\
\noindent \textbf{Proof}. See \cite{Hoy18}, Theorem 15.\\\\
\textbf{Diamond reformulation}. The rough idea is the following:\\\\
Conjecture \cite{Dob21a}. Let $S$ be a perfectoid space, $\mathcal{D}^{\diamond}$ a stable dualizable presentable
category, and $R$ a p-adic shtuka on $X$. Let $\mathcal{D}$ be a stable compactly generated category and $F $: Cat$^{idem}_{St} \rightarrow \mathcal{T}$ a localizing invariant that preserves filtered colimits. Then:\\

$F_{cont}(Sht(\mathbb{S}^n, \mathcal{C})) \simeq \Omega^nF_{cont}(\mathcal{D}^{\diamond})$ with $S$ $SpaQ_p \times SpaQ_p$. 


\section{Implications}
\label{sec:consequences}

\subsection{A language scheme}
\label{subsec:lang}
Presently, we do not have a way to talk (or write) about $n$-awareness -- our language is antiphrastical. The Danish philosopher Søren Kierkegaard once remarked that piling up (empirical) knowledge in hope to get a glimpse of the nature of being commits one to ``a perpetually self-repeating false sorites'' \cite{Kie38}. One could interpret this as the problem of arriving at knowledge about subjective experience by way of chaining together (empirical) propositions -- a version of the explanatory gap argument \cite{Lev83}. An alternative, following Kierkegaard, would be to repeatedly live through one's subjective experience, thereby expressing its nature.

Ontological statements could thereby be reconceived as being linguistic statements.  We are ``upgrading'' Wittgenstein’s language games \cite{Wit53} to a ``language scheme''. 
This language scheme follows from our meta-model. Instead of saying that meaning (reference) is a derivative from whatever game is at play, we say that meaning (reference) is a derivative from whatever scheme is at play. In view of Conjecture 1, the language scheme refers to a ``profinite language'' with expressions that consist of profinitely many words that are ``glued together''.%
\footnote{
What is the glossematics of such a profinite language? Recall that a profinite set $S$ is called extremely disconnected if the closure of any open subset $U \rightarrow S$ is still open.
The philosopher and linguist Roland Barthes stated that ``writing ceaselessly posits meaning ceaselessly to evaporate it, carrying out a systematic exemption of meaning.''  \cite{Bar77}
So, where is the meaning?  We propose that word-meanings are a delicate thing: while we have profinitely many copies of words that are ``glued together'', they never refer unambiguously to a single meaning ``out there''. The glossemes correspond to the geometric impurities of the perfectoid diamond, which is the glossematics.
}


%
A basic example for language misuse is the stipulation of a substantive ``self'' that endures through ``time'', which might in actuality merely reflect the grammatical (and indexical) structure of our ordinary language. (Even though such a ``self'' does not refer to an existing ``thing out there'', it is highly useful, e.g., for communicating certain ideas, desires or intentions.) Another example would pertain to the idea of some ``private internal language''. Actual languages depend on them being used by a community.\footnote{This is also true for ``biological languages'' implemented in a living ``community'' making up a body.} 
The meaning of a word does not refer to an essentially private ``thing''.
Analogously, one could proceed in the case of a language scheme. 

First, we outline the basic grammatical structure of this scheme. 
Modern English language contains three grammatical cases (subjective, objective, and possessive) with different declensions for each case. We extend this grammatical structure by adding a number to the respective pronoun, indicating the simultaneous presence of different ``me's'' within an experience of multiple awareness (Tab. \ref{tab1}). Formally, the resulting “$n$-declension” allows for all combinations of “$n$, $n-1$, $\hdots$, 1-declensions” to be present in a sentence. So for example, a statement expressing double-awareness could contain a “1-I” as well as a superposed “2-I” as subjects. Analogously,  “$n$-conjugation” could be defined. 

Second, we envision a novel model of temporality. Spatio-temporally multiplied awareness is obviously not reconcilable with theories which posit awareness to be bound to one particular location in  (physical) space and time. However, non-physicalist theories (outside of space-time) seem to be able to accommodate the concept of $n$-awareness. Many mystics in the religious traditions from the East and West reported a similar kind of experience. If awareness is \textit{not} bound to a single region in space and time, this also suggests that awareness cannot be understood as an emergent property of localized physical systems.

Our model conceives of a temporal multiplicity with a categorified model of $n$-time, evidenced in the proposed language scheme by way of $n$-inflection for $n$-conjugation. Modern English is spoken in local, linear time, yet it allows the inflectional change of verbs by way of conjugation. We extend the idea of language spoken in linear time, conjugated over three tenses, to one of $n$-conjugation as follows: Instead of using only past, present, future, and their perfect correspondences, present-perfect, past-perfect, future-perfect, we allow for 1-past, 2-past, $\hdots$, 1-present, 2-present, $\hdots$, 2-present-perfect, 2-future-perfect, $n$-future-perfect etc, which is what we call $n$-time. This generalizes the discussion on temporal experience, which has traditionally been expressed in terms of a ``tensed'' experienced time vis-a-vis an ``untensed'' physical parameter time \cite{Tag08}. 
\begin{sidewaystable}
\centering
\caption{Grammar of a language scheme for $n$-awareness, extended from ordinary English case structure.}
\label{tab1}
\begin{tabular}{llllllll}
\\\toprule
&  & & $n=1$ & $n=2$ & $n=3$ & $ \cdots$ &$n=k$\\\hline
\textbf{subjective} & I (we) & $\rightarrow$  & 1-I (1-we) & 1-I (1-we) & 1-I (1-we) & $\cdots$ & 1-I (1-we) \\
& & & &  2-I (2-we) & 2-I (2-we) && 2-I (2-we) \\
& & & &  & 3-I (3-we) && $\vdots$ \\
& & & & & & & $k$-I ($k$-we)\\\\

\textbf{objective }& me (us) & $\rightarrow$  & 1-me (1-us) & 1-me (1-us) & 1-me (1-us) & $\cdots$ & 1-me (1-us) \\
& & & &  2-me (2-us) & 2-me (2-us) && 2-me (2-us) \\
& & & &  & 3-me (3-us) && $\vdots$ \\
& & & & & & & $k$-me ($k$-us)\\ \\

\textbf{dep. possessive} & my (our) & $\rightarrow$  & 1-my (1-our) & 1-my (1-our) & 1-my (1-our) & $\cdots$ & 1-my (1-our) \\
& & & &  2-my (2-our) & 2-my (2-our) && 2-my (2-our) \\
& & & &  & 3-my (3-our) && $\vdots$ \\
& & & & & & & $k$-my ($k$-our)\\ \\

\textbf{indep. possessive} & mine (ours) &$\rightarrow$  & 1-mine (1-ours) & 1-mine (1-ours) & 1-mine (1-ours) & $\cdots$ & 1-mine (1-ours) \\
& & & &  2-mine (2-ours) & 2-mine (2-ours) && 2-mine (2-ours) \\
& & & &  & 3-mine (3-ours) && $\vdots$ \\
& & & & & & & $k$-mine ($k$-ours)\\ \\

\end{tabular}
\end{sidewaystable}
\subsection{Revisiting problems of temporality}
\label{sec:time}
Different temporal ontologies have been proposed throughout the ages, without coming to definite conclusion. Rather than proposing yet another metaphysical framework, we wish to concentrate on particular problems that feature prominently in recent and not so recent discussions. 
We hope that approaching these problems structurally, inspired by our meta-model, will lead to a remedy.
The problems we wish to look at are the following:
\begin{enumerate}
\item \textbf{How is change possible?}
Change manifests itself every day. But what does it refer to exactly? The idea that ``change'' does not really exist goes back to the works of Parmenides and the Eleatic school of philosophy. While this philosophy has been influential up to this day, for example, in the metaphysical thought of Martin Heidegger  \cite{Hei42}, our everyday experience seems much better captured in the ``everything flows'' of Heraclitus \cite{Die52}.

The American philosopher David Lewis revived the problem of change for the philosophy of time in the problem of temporary intrinsics\footnote{An intrinsic property is a property an object has irrespective to the relations it bears to other objects. The converse would be a relational property: being an uncle is not a property that I have independent of my nephews. } \cite{Lew86}:
\begin{quote} 
For instance shape: when I sit, I have a bent shape; when I stand, I have a straightened shape. Both shapes are temporary intrinsic properties; I have them only some of the time. How is such change possible? 
\end{quote}
According to David Lewis the problem of temporary intrinsics could be solved in three ways. Either one acknowledges that there \textit{are no} intrinsic properties, just ``disguised relations''; or one believes that only those properties that exist at the \textit{present} moment are real, whereas the properties that an object seem to have had previously are, in some sense, fictional (this position is known as ``presentism''); or one accepts that objects have genuine temporal parts (e.g. the me-yesterday, the me-now, and the me-tomorrow). The latter solution to the problem of temporary intrinsics has been deemed the only viable solution to the problem of temporary intrinsics which is not ``incredible'' \cite{Lew86} and started the appreciation of ``perdurance'' theories in the modern philosophy of time at the expense of so-called ``endurance'' theories that conceive of persisting wholes without temporal parts \cite{Haw20}. In addition to Lewis' metaphysical rejunevation, much support for a ``perdurance-like'' theory seems to come from science, in particular Einstein's theory of relativity. Perdurance theory, so it is often but not always believed, squares well with the belief that space-time forms a four-dimensional continuum as described by the special theory of relativity \cite{Haw20}.
\item \textbf{What is simultaneity?}
The second problem worth mentioning in this respect is the problem of simultaneity (coincidence). It seems that, when discussing $n$-awareness we postulate the simultaneous presence of two experiences. This seems to violate the basic intuition that no two objects could occupy the same place in time unless they are the same (or unless they share temporal parts: the statue and the clay have temporal parts that overlap\footnote{Another argument in favor of perdurantism.}). But it also seems to be in conflict with basic principles of physics according to which there can be no ``absolute'' notion of simultaneity. 

\item \textbf{Is synchronous reference possible?}
While perdurantism claims to solve problems of temporal coincidence (which are only problems if designed to be so), it has its own problems when trying to account for the acts of synchronous and asynchronous reference. For instance, if I touch my hand on a hot stove and get burned, it takes ``time'' for my system to register ``pain''.  My feeling pain ``now'' is a result of my action in the immediate past. One could claim that all present feelings are results of actions in the past, be it immediate or not. As such, there is no such thing as synchronous reference, only asynchronous. This means it is important to question how received ontologies (such as perdurantism or endurantism) handle the question of asynchronous reference. At least prima facie, it seems they cannot, as observed by Bertrand Russell \cite{Rus21}: 
\begin{quote}
There is no logically necessary connection between events at different times; therefore nothing that is happening now or will happen in the future can disprove the hypothesis that the world began five minutes ago.
\end{quote}





\end{enumerate}
%
We replace the idea that objects are bearers of (intrinsic or extrinsic) properties by the idea that ``properties'' are identical to the morphisms (relations) within a category. ``Change'' in this context is tantamount to the addition (removal) of morphisms.\\
The traditional distinction between perdurantism and endurantism could be illustrated by the notion of a so-called ``coequalizer'':
\begin{equation}
\nonumber
\xymatrixrowsep{1cm}
\xymatrixcolsep{1.5cm}
\xymatrix{
X \ar@<1ex>[r]^f \ar@<-1ex>[r]_g & Y \ar[r]^q \ar[rd]_{q'} & Q \ar@{..>}[d]^u \\
& & Q'
}
\end{equation}
In category theory, a coequailizer refers to a single object (a ``colimit'') associated to the different morphisms $f$ and $g$ between objects $X$ and $Y$, such that $q \circ f = q\circ g$. Furthermore, the objects $Q$ is universal, meaning it is unique ``up to'' an isomorphism $u$.
It follows that \textit{properties} (i.e. morphisms within a category) are associated to a single and unique (up to isomorphism) \textit{object}. Whereas perdurantism, translated into the language of category theory, is about change between such properties (i.e. the addition of new morphisms), endurantists refer to the unchanging (persisting) object defined by them. 
The diamond $SpdQ_p = Spa(Q_p^{cycl})/ \underline{Z_p^\times} $ is the coequalizer of $\underline{Z_p^\times} \times Spa(Q_p^{cycl})^\flat \rightrightarrows Spa(Q_p^{cycl})^\flat$,  where one map is the projection and the other is the action; cf. \cite{SW20}.

%
Our solution 
superficially seems to correspond to a perdurantist representation with the $n$-declension of $1$-her, $2$-hers, ... $n$-her, but these are relational properties that do not necessarily refer to temporal or intrinsic parts.  It can be asked how different Mary-tomorrow is from Mary-today given that Mary-tomorrow has more morphisms?  If Mary refers to the objects of a category and the properties are its morphisms, then saying that Mary has ``changed'' is merely to say that Mary has added connections/morphisms. Else, we can say that Mary-today is the same as Mary tomorrow ``up to'' isomorphism.

Analogous to how perdurantists resolve the problem of coincidence by noting that temporal parts can indeed ``overlap'' without implying that the two objects that overlap are identical, we note that categories (given they have at least some basic structural features: e.g. they are topological) could too be said to ``overlap''. But this does not commit us to perdurantism as an ontological position.

Note that the problem of simultaneity is mainly a (conceptual) ``design issue'' that stems from a linear notion of time, where simultaneity is conceived in terms of (``temporal'') coincidence, or, alternatively, from the treatment of time in the framework of Minkowski space-time (often called a ``fourth dimension''). We instead choose to model time in terms of an equivalence relation using homotopy theory -- ``time'' is not an absolute (ontological) notion, but instead refers to a relative (epistemic) ``ordering'' scheme of experiences. Thus, to every level of awareness, there corresponds a level of time.  $1$-awareness corresponds to $1$-time; $2$-awareness corresponds to $2$-time; etc. 
\\
We offer a structural notion of ``simultaneity of experience'' by identifying it with the commutativity of diagrams in a homotopy category, and respectively by the ``up to'' notion. Commutativity classifies the equivalence of all possible ways to get to a destination.\footnote{Take, for example, the chain map between complexes: Fig. \ref{fig:chain}, in the appendix.}
There is no (structural) difference in choosing one way over the other. There is no indicated starting point or canonical progression. Rather, all possible paths are ``revealed at once'' (and even infinite paths might be alluded to). The ``up to'' notion grants a relative notion of equivalence that corresponds to commutativity. Making the statement that all morphisms are equivalent ``up to'' homotopy means that they are equivalent with respect to homotopy. There is no substantial way to distinguish one morphism over any other. In a sense, commutativity is our way of geometrizing the ``up to'' notion. Equivalence and hence simultaneity is never truly absolute. Simultaneity, on our view, expresses a homotopy equivalence and homotopy equivalences are neither perdurantist nor endurantist. 

Moreover, our meta-model provides a structural framework for (a)synchronous reference.
%
That there is only asynchronicity is implicit in the notion of diamond nonlocality \cite{Dob21b}. 
Diamond nonlocality is a perfectoid version of nonlocality that arises from the nontrivial geometry of the profinitely many copies of $Spa(\mathcal{C})$. 
The asynchronicity comes from the profinite condition. There is no canonical reference frame because of the existence of many such quasi-pro-\'etale covers, the pullback of which gives the profinitely many copies. Events and their seeming synchrocity lie entirely in the quasi-pro-\'etale covers.\footnote{Due to this mathematically ``holographic" structure of the diamond, the idea is to reconstruct the AdS/CFT holographic principle \cite{juan} using diamonds, and the six operations in Scholze's \'etale cohomology of diamonds \cite{Sch17}. The diamonds represent the conformal field theory. The six operations reconstruct an analytic version of what is Anti de Sitter space.}

%

\subsection{Self-reflexive selves}
\label{sec:selves}

Two properties typically associated with the notion of self are, first, that the conscious experience of a self is fundamentally self-reflexive: each moment of an experience always refers to the experiential context as a whole: Experience is embedded into a self that has these experiences. If a self is nothing but a ``bundle''  of such experiences (or more technically: a weakly persistent moduli space \cite{msp} of relations between them), any moment of this experience refers to the whole bundle. This is not necessarily the case when we talk about ``minimal experiences as such'' \cite{Met20}, but it should be the case when experience is associated with a single self that has them.%
\footnote{Whether there truly is such an ``experience as such'' outside the context of a self, is still an open debate.}
Second, we are interested in a model that acknowledges plurality and individuality, something that seems at odds with the idea of a metaphysical monism that postulates all things to ``be of the same kind'' (at best, individual selves were instances of a universal natural process).

%

If we work in an $(\infty,1)$-topos, which is presentable and satisfies a descent condition, then it can be said that experience is ``self-reflexive'' with respect to the objects it contains, starting from the original representations (i.e. the moments of 1-awareness). 
%
Phenomenologically, this could be taken to mean that my awareness now is (in some sense) equivalent to my awareness yesterday, and even worse: \textit{your} awareness tomorrow is (in some sense) equivalent to \textit{my} awareness two days ago. This conflicts with our intuition that our experiences (across individuals but also across times) are \textit{unique}. The notion of equivalence ``up to'' homotopy affords to deflate this strong notion of equivalence. 

You are not strictly the same ``you'', with some substantialist notion of self in the background. But, using our model of homotopy types, you could be homotopy equivalent to a ``structural equivalent'' of yourself.
In our meta-model, the structural equivalence lies in taking the Efimov K-theory of the localization sequence of the diamond representing awareness to get the global view of equivalence.%
\footnote{Informally, this structural equivalence uses what one of the authors calls ``perfectoid entanglement". ``Perfectoid entanglement entropy'' is proposed as a measure of the degree of equivalence  (cf. \cite{Dob21b}).}


Considering the results of subsection \ref{sec:time}, there is no conception of a self ``enduring over time'', as if time were somehow exterior to experience; as if time were a constantly present  fluid through which selves move (however \textit{that} happens); or as if time were a simple counter of experiences. 

\section{Summary}

In section \ref{sec:model}, we outlined a framework based on categories of representations and looked at hierarchies that arise from the category theoretic treatment, starting with the morphisms defining ``$1$-awareness'' up to (invertible) ```$n$-morphisms''. We emphasized the topological nature of the encountered relations between moments of awareness, and we postulated to model this with homotopy theory and the theory of $(\infty,1)$-topoi. This was our basic (``vanilla'') framework for studying subjective experience. 

We then enriched the framework with three conjectures that define our meta-model. We propose a concrete mathematical structure, the ``perfectoid diamond'', as a model of subjective experience and appeal to Efimov K-theory to study the equivalence relations across all such diamonds. Specifically:

\begin{itemize}

\item \textbf{Conjecture 1}. The  $(\infty,1)$-category of perfectoid diamonds is an $(\infty,1)$-topos.  

\item \textbf{Conjecture 2}. Topological localization, in the sense of Grothendieck-Rezk-Lurie $(\infty,1)$-topoi, extends to the $(\infty, 1)$-category of diamonds.

\item \textbf{Conjecture 3}.  The meta-model takes the form of Efimov K-theory of the $(\infty,1)$-category of diamonds.
\end{itemize}

Finally, some of the implications of our meta-model were discussed in the ensuing section \ref{sec:consequences}, where we derived the grammar for a novel ``language scheme'' (including the notions of ``$n$-declension'' and ``$n$-time''), revisited some traditional problems in the philosophy of time, and eventually investigated a notion of a self as a ``bundle'' of experiences.
These discussions had the purpose of, first, illustrating the framework in terms of language and generalizing a point made previously by various philosophers: how the intricacies of our language wrongly suggest to assume substantive entities (a disembodied ``self''; ``meanings'' out there), where there in fact are none. With reference to the perfectoid diamond model, we proposed a ``pro-finite'' language with expressions that consists of profinitely many words that are ``glued together'' without ever carrying some unambigious reference.

Second, and related, many problems in the philosophy of time are artifacts of such a confused language use. The remedy to these problems might not be found in a metaphysical solution (e.g. ``perdurantism''), but in a re-conceptualization of problematic terms such as ``simultaneity''. We moreover hinted at a natural resolution of the problem of (a)synchronous reference within our meta-model.

 Third, the notion of a ``self'' could be understood in terms of a bundle of experiences that is ``self-reflexive" and possess a form of ``individuality''. The notion of ``structural equivalences" between classes of experiences is made precise in the Efimov K-theory of perfectoid diamonds.
Note that in each case certain supposedly ``real'' entities (linguistic meanings, change, selves) were reduced to structuralist notions.

We consider our model to reflect a new movement in mathematics which seeks to make functional analysis a branch of commutative algebra and therefore provide a geometrization of analytical techniques. We conceive of replacing our model's foundation of topological spaces with that of condensed mathematics \cite{con}. The study of subjective experience is proposed to be a geometrical, rather than a computational, project. 
%
%
%
%
%

\clearpage

\subsection*{Appendix}
\label{appendix}

An important invariant of a mathematical space is encoded by its \textit{homology} group - the number of holes in that space. As such, they provide a means to compare spaces. For $X$ a topological space, a set of topological invariants $H_0(X)$, $H_1(X)$,..., called the homology groups of $X$, represent the homology of $X$. The number of $k$-dimensional holes in $X$ is encoded by the $k$th Homology group $H_k(X)$. For instance, $H_0(X)$ encodes the ``path connected'' components of $X$, where a (0-dimensional) hole encodes if the space is disconnected. 

As an example, let us examine the homology groups of $S^1$, the 1-dimensional sphere (which is really just a circle). Take $X$ to be $S^1$. $X$ is connected and has one 1-dimensional hole and no other holes for $k>1$. The homology groups of $X$ take the form:  
\begin{equation}
H_{k}(S^{1})={\begin{cases}\mathbb {Z} &k=0,1\\\{0\}&{\text{otherwise}}\end{cases}}
\end{equation}
Take $X$ to be $S^2$, the 2-dimensional sphere (which is just the surface of a ball). $S^2$ is connected and has just one 2-dimensional hole.  The homology groups of $X$ are represented as:  
\begin{equation}
 H_{k}(S^{2})={\begin{cases}\mathbb {Z} &k=0,2\\\{0\}&{\text{otherwise}}\end{cases}}
 \end{equation}
 \\
 To introduce the derived categorical setting, we will first explain the spirit of \textit{quasi-isomorphisms} using the notion of quasi-categories formulated by Andre Joyal \cite{Joy08}. A quasi-isomorphism is a morphism $A \rightarrow B$ of chain complexes\footnote{Respectively of cochain complexes.} such that the induced morphisms $H_n(A, \cdot) \rightarrow H_n(B,\cdot), H^n (A,\cdot) \rightarrow H^n(B, .)$ of homology groups\footnote{Respectively of cohomology groups.} are isomorphisms for all $n$. 

Quasi categories are homotopoi \cite{Joy08} which possess rich general structures and do not necessarily have a uniquely defined composition of morphisms. Quasi-categories are like ordinary categories in that they are certain simplicial sets which contain objects (the 0-simplices of the simplicial set) and morphisms between these objects (1-simplices). Unlike categories, however, morphisms can be composed, but the composition is well-defined only up to still higher order invertible morphisms. This means that all possible morphisms which serve as the composition of two 1-morphisms are related to each other by 2-morphisms called 2-simplices, which resemble homotopies. It turns out that every Kan Complex \cite{kan} is a quasi-category. So from this Kan Complex we build the notion of the chain complex, which we will use in the next section.

We mention two reasons that we work in the \textit{derived categorical }setting. 
One reason is that knowing the homology of a space does not give complete information about its homotopy type. This is seen by that fact that there exist topological spaces $X$ and $Y$ such that $H_i(X)$ is isomorphic to $H_i(Y)$ for every $i$, but $X$ is not homotopy equivalent to $Y$.  The derived category remembers the entire complex, which is crucial to our model of $n$-awareness, and the consequent model structure gives us nice classes of morphisms which axiomatize homotopy theory. 
Another reason is that the derived category setting allows us to \textit{localize} in the category setting. Localization is a formal process of adding inverses to a space. A category can be localized by formally inverting certain morphisms, such as the weak equivalences in the homotopy category of a model category. We use a special case of localization called Bousfield localization \cite{bou}, which assigns a new model category structure with more weak equivalences to a given model category structure. So Bousfield localization allows us to make more morphisms count as weak equivalences and this is formally how we get from 1-awareness to n-awareness. 

\begin{figure}
\begin{equation}
\nonumber
\xymatrixrowsep{1.5cm}
\xymatrixcolsep{1.5cm}
\xymatrix{
X \ar[r]^{\tilde{f}_0} \ar[d]_{X \times \{0\}} & E \ar[d]^{\pi}  & Y & \ar[l]_{\tilde{f}_0} X\\
X \times I \ar[r]^{f} \ar@{.>}[ur]^{\tilde{f}} & B & Y^I \ar@{->>}[u]^{p_0} & A \ar[u]_{i} \ar[l]_{f_0} \ar@{.>}[ul]_{\tilde{f}}
}
\end{equation}
\caption{Homotopy lifting and extension property for topological spaces $E$ and $B$. In the leftmost figure, the homotopy lifting property allows homotopies in the space $B$ to be uplifted to the space $E$ for any homotopy $f: X \times [0,1] \rightarrow B$ and for any map $\tilde{f}_{0}: X \rightarrow E$ such that $f_0 = \pi \circ \tilde{f}_{0}$. A lifting $\tilde{f}$ corresponds to a dotted arrow giving a commutative diagram. In the rightmost figure, the homotopy extension property extends certain homotopies defined on a subspace to a larger space. The homotopy extension property is dual to the homotopy lifting property \cite{hom}. 
}  
\label{fig:lifting}
\end{figure}
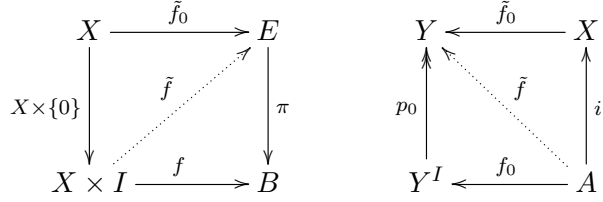

To axiomatize homotopy theory, we use the construction of a Quillen model structure \cite{mod}. A model structure on a category consists of three classes of morphisms: weak equivalences, fibrations, and cofibrations. Weak equivalences are quasi-isomorphisms, maps which induce isomorphisms in homology. Cofibrations are maps that are monomorphisms that satisfy the homotopy extension property. Fibrations are maps that are epimorphisms that satisfy homotopy lifting property (Fig \ref{fig:lifting}). In the derived setting, quasi isomorphisms are used as the class of weak equivalences, fibrations mimic surjections, and the cofibrations mimic inclusions. From this model structure, we will define the notion of simultaneity. 

We let $A$ be a Grothendieck abelian category, such as the category of abelian groups.\footnote{Grothendieck worked on unifying various constructions in mathematics. For instance, the Grothendieck group construction is the most universal way of constructing an abelian group from a commutative monoid.} The Grothendieck abelian category is an $AB_5$ category with a generator. $AB_5$ categories are $AB_3$ categories (abelian categories possessing arbitrary coproducts) in which filtered colimits of exact sequences are exact \cite{Gro60}.
The category of abelian groups is a prototypical example of a Grothendieck category, with generator the abelian group $\mathbb{Z}$ of integers. The category of abelian groups has as objects abelian groups and as morphisms group homomorphisms. We use Grothendieck categories because we need a category universally enriched over abelian groups to model n-awareness. 
Groups encode symmetries. But what is 1-awareness an symmetry of? We hypothesize that it encodes a complexity class of Turing degree 0 \cite{tur}.

We then construct a new ``derived category'', $D(A)$, whose objects are complexes of objects of $A$ and whose morphisms are chain maps. $D(A)$ contains a model structure that will be our model of $n$-awareness.\footnote{For a more detailed exposition see the work of A. Caldararu \cite{Cal05}.} 

Firstly, we define a chain complex. A chain complex $(A_\bullet, d_\bullet)$ is a sequence of abelian groups  $..., A_0, A_1, A_2, A_3, A_4, ...$ connected by homomorphisms (called boundary operators or differentials) $d_n: A_n \rightarrow A_{n-1}$, such that the composition of any two consecutive maps is the zero map. Explicitly, the differentials satisfy $d_n \circ d_{n+1} = 0$, or with indices suppressed, $d^2 = 0$. A chain complex has the form:
\begin{equation}
\nonumber
{\displaystyle \cdots {\xleftarrow {d_{0}}}A_{0}{\xleftarrow {d_{1}}}A_{1}{\xleftarrow {d_{2}}}A_{2}{\xleftarrow {d_{3}}}A_{3}{\xleftarrow {d_{4}}}A_{4}{\xleftarrow {d_{5}}}\cdots }
\end{equation}
Secondarily, a chain map $f$ between two chain complexes $(A_\bullet, d_\bullet), (B_\bullet, d_\bullet)$ is a sequence $f_n$ of homomorphisms $f_n: A_n \rightarrow B_n$ for each $n$ that commutes with the differentials on the two chain complexes, where  $d_B, n \circ f_n = d_A,n \circ f_{n-1}$. A chain map takes the form of the commutative diagram in Fig. \ref{fig:chain}.
$(f_{\bullet })_{*}:H_{\bullet }(A_{\bullet },d_{A,\bullet })\rightarrow H_{\bullet }(B_{\bullet },d_{B,\bullet })$ on preserves cycles and boundaries, so $f$ induces a map on homology.   
\begin{figure}[htb]
\begin{equation}
\nonumber
\xymatrixrowsep{1.25cm}
\xymatrixcolsep{1.25cm}
\xymatrix{
\hdots & A_{n-1} \ar[l]_{d_{A,n-1}} \ar[d]_{f_{n-1}} & A_{n} \ar[l]_{d_{A,n}} \ar[d]_{f_{n}}  & A_{n+1} \ar[l]_{d_{A,n+1}} \ar[d]_{f_{n+1}} & \hdots \ar[l]_{d_{A,n+2}} \\
\hdots & B_{n-1} \ar[l]_{d_{B,n-1}} & B_{n} \ar[l]_{d_{B,n}}  & B_{n+1} \ar[l]_{d_{B,n+1}} & \hdots \ar[l]_{d_{B,n+2}}
}
\end{equation}
\caption{A commutative diagram of a chain map between two chain complexes $(A_\bullet, d_\bullet), (B_\bullet, d_\bullet)$, with a sequnce of morphisms $f_n: A_n \rightarrow B_n$ and differentials $d_B, n \circ f_n = d_A,n \circ f_{n-1}$. }
\label{fig:chain}
\end{figure}


\clearpage

\begin{thebibliography}{}



\bibitem{Dob21a} Dobson, S. \textit{Efimov K-theory of Diamonds}, in preparation.

\bibitem{kth} ncatlab authors. K-theory. Available online: \url{https://ncatlab.org/nlab/show/K-theory} (accessed February 10, 2021).

\bibitem{SW20} Scholze, P and Weinstein, J. \textit{Berkeley Lectures on P-adic Geometry}, Princeton University Press, Annals of Mathematics Studies
Number 207.

\bibitem{Hoy18} Hoyois, M., \textit{K-theory of Dualizable Categories (after rA. Efimov)}, notes http://www.mathematik.ur.de/hoyois/papers/efimov.pdf {accessed January 11, 2021}.

\bibitem{Dob21b} Dobson, S., \textit{Perfectoid Quantum Physics and Diamond Nonlocality}, in preparation.

\bibitem{Tye09} Tye, M. "Representational Theories of Consciousness", in \textit{The Oxford Handboook of Philospohy of Mind}, edited byMcLaughlin, B.P., Beckermann, A., and Walter, S., Oxford University Press, 2009, pp. 253--267.
\bibitem{Hei27} Heidegger, M. \textit{Sein und Zeit}, Niemeyer, 1927.

\bibitem{Spi77} Spinoza, B. \textit{Ethica, ordine geometrico demonstrata}, 1677.

\bibitem{per} ncatlab authors. Perfectoid Space. Available online: \url{https://ncatlab.org/nlab/show/perfectoid+space} (accessed January 29, 2021).

\bibitem{hub} ncatlab authors. Huber Space. Available online: \url{https://ncatlab.org/nlab/show/Huber+space} (accessed January 29, 2021).

\bibitem{pro} ncatlab authors. Profinite Space. Available online: \url{https://ncatlab.org/nlab/show/profinite+space} (accessed January 29, 2021).

\bibitem{geo} ncatlab authors. Geometric Point. Available online: \url{https://ncatlab.org/nlab/show/geometric+point} (accessed January 29, 2021).


\bibitem{gro} ncatlab authors. Grothendieck Topology. Available online: \url{https://ncatlab.org/nlab/show/Grothendieck+topology} (accessed January 29, 2021).

\bibitem{Far16} Fargues, L., \textit{Geometrization of Local Correspondence, an Overview}, arXiv:1602.00999 [math.NT] (accesseed February 10, 2021).

\bibitem{Wit53} Wittgenstein, L. \textit{Philosophical Investigations}, Macmillan Publishing Company, 1953.

\bibitem{Hug18} Huggett, N, and Hoefer, C. ``Absolute and Relational Theories of Space and Motion'', in: \textit{The Stanford Encyclopedia of Philosophy (Spring 2018 Edition)}, edited by Zalta, EN. \url{https://plato.stanford.edu/archives/spr2018/entries/spacetime-theories/} (accessed January 29, 2021).

\bibitem{Lei14} Leibniz, GW. \textit{La Monadologie}, 1714.

\bibitem{upt} Wikipedia authors. Up to. Available online: \url{https://en.wikipedia.org/wiki/Up_to} (accessed January 29, 2021).

\bibitem{Hum40} Hume, D. \textit{A Treatise of Human Nature}, 1740.

\bibitem{abe} Wikipedia authors. Abelian group. Available online: \url{https://en.wikipedia.org/wiki/Abelian_group} (accessed January 29, 2021).

\bibitem{com} ncatlab authors. Complex. Available online: \url{https://ncatlab.org/nlab/show/complex} (accessed January 29, 2021).

\bibitem{Art63} Artin, M, Grothendieck, A, and Verdier, J-L. \textit{Th\'eorie des topos et cohomologie \'etale des sch\'emas, S\'eminaire de G\'eom\'etrie Alg\'ebraique du Bois-Marie 1963-64}, dirig\'e par M. Artin, A. Grothendieck, J.-L. Verdier, SLN 269, 270, 305, Springer-Verlag, 1972, 1973.









\bibitem{Ehr15} Ehresmann AC, and Gomez-Ramirez. J. Conciliating neuroscience and phenomenology via category theory. \textit{Prog. Biophys.} \textbf{2015}, \textit{119}, 347-359.

\bibitem{Tsu16} Tsuchiya, N, Taguchi, S, and Saigo, H. Using category theory to access the relation between consciousness and the integrated information theory, \textit{Neuroscience Research} \textbf{2016}, 107, 1–7. 

\bibitem{Pre19} Prentner, R. Consciousness and Topologically Structured Phenomenal Spaces, \textit{Consciousness \& Cognition} \textbf{2019},  \textit{70}, 25--38.

\bibitem{Nor19} Northoff, G, Tsuchiya, N, and Saigo, H. Mathematics and the Brain: A Category Theoretical Approach to Go Beyond the Neural Correlates of Consciousness. \textit{Entropy} \textbf{2019}, 21(12), 1234–21. 

\bibitem{Kle20} Kleiner, J. Mathematical Models of Consciousness, \textit{Entropy} \textbf{2020}, \textit{22(6)}, 609.

\bibitem{Sig20} Signorelli CM, Wang, Q, and Kahn, I. A Compositional Model of Consciousness based on Consciousness-Only, 2020. \url{https://arxiv.org/abs/2007.16138} (accessed January 29, 2021).

\bibitem{Ehr07} Ehresmann AC, and Vanbremeersch, J-P. \textit{Memory Evolutive Systems: Emergence, Hierarchy, Cognition.} Elsevier, 2007.

\bibitem{Put67} Putnam, H. ``Psychological Predicates'',  in \textit{Art, Mind, and Religion}, edited by Capitan, WH, and Merrill, DD. University of Pittsburgh Press, 1967, pp. 37--48.

\bibitem{Hus75} Husserl E. \textit{Logische Untersuchungen}, Husserliana XIII, 1975.

\bibitem{Blo96} Block, N. What is functionalism? \textit{The Encyclopedia of Philosophy Supplement}, 1996. \url{http://www.nyu.edu/gsas/dept/philo/faculty/block/papers/functionalism.pdf} (accessed Jan 26, 2021).

\bibitem{Wu18} Wu, W. "The Neuroscience of Consciousness", in \textit{The Stanford Encyclopedia of Philosophy (Winter 2018 Edition)}, edited by Zalta, EN. \url{https://plato.stanford.edu/archives/win2018/entries/consciousness-neuroscience/}  (accessed January 29, 2021).

\bibitem{Cra13} Crane, T. \textit{The Objects of Thought}, Oxford University Press, 2013.

\bibitem{Hof08} Hoffman, DD. Conscious realism and the mind-body problem. \textit{Mind and Matter} \textbf{2008}, \textit{6(1)}, 87–121.

\bibitem{Hof14} Hoffman, DD, and Prakash, C. Objects of Consciousness, \textit{Frontiers in Psychology} \textbf{2014}, 5:577.

\bibitem{Gof17} Goff, P. \textit{Consciousness and Fundamental Reality}, Oxford University Press, 2017.

\bibitem{Luc04} Illusie, L. What is a topos? \textit{Notices of the AMS} \textbf{2004}, Vol 51, Number 5.

\bibitem{top} ncatlab authors. Topos. \url{https://ncatlab.org/nlab/show/topos} (accessed January 29, 2021).

\bibitem{triangulated} ncatlab authors. Triangulated Category. \url{https://ncatlab.org/nlab/show/triangulated+category} (accessed January 29, 2021).

\bibitem{hig} ncatlab authors. Higher Category Theory. \url{https://ncatlab.org/nlab/show/higher+category+theory} (accessed January 29, 2021).

\bibitem{der} ncatlab authors. Derived Category. \url{https://ncatlab.org/nlab/show/derived+category} (accessed January 29, 2021).

\bibitem{groth} ncatlab authors. Grothendieck Topology. \url{https://ncatlab.org/nlab/show/Grothendieck+topology} (accessed January 29, 2021).

\bibitem{inf} ncatlab authors. (infinity,1)-topos. \url{https://ncatlab.org/nlab/show/(infinity,1)-topos} (accessed January 29, 2021).

\bibitem{Kie38} Kierkegaard, S. Journal entry from August 17, 1838.

\bibitem{Lev83} Levine, J. Materialism and Qualia: The Explanatory Gap. \textit{Pacific Philosophical Quarterly} \textbf{1981}, 64, 354–361.

\bibitem{dec} ncatlab authors. Descent. \url{https://ncatlab.org/nlab/show/descent} (accessed January 29, 2021).

\bibitem{Bar77} Barthes, R. The Death of the Author, in: \textit{Image, Music, Text}, Fontana Press, 1977.

\bibitem{Tag08} McTaggart, JE. The Unreality of Time, \textit{Mind}, 1908.

\bibitem{Hei42} Heidegger, M. \textit{Parmenides, lecture in the spring semester 1942/43}, Vittorio Klostermann, 2018.

\bibitem{Die52} Diels, H, and Kranz, W. \textit{Die Fragmente der Vorsokratiker}, Weidmann, 1952.

\bibitem{Lew86} Lewis, DK. \textit{On the Plurality of Worlds}, Blackwell, 1986.

\bibitem{Haw20} Hawley, K. ``Temporal Parts'', in \textit{The Stanford Encyclopedia of Philosophy (Summer 2020 Edition)}, edited by Zalta, EN. \url{https://plato.stanford.edu/archives/sum2020/entries/temporal-parts/}. (accessed January 29, 2021)

\bibitem{Rus21} Russell, B. \textit{The Analysis of Mind}, G. Allen \& Unwin Limited: London, 1921.

\bibitem{msp} Wikipedia authors. Moduli space. \url{https://en.wikipedia.org/wiki/Moduli_space} (last accessed January 29, 2021)

\bibitem{Met20} Metzinger, T. Minimal phenomenal experience. Meditation, tonic alertness, and the phenomenology of ``pure” consciousness. \textit{Philosophy and the Mind Sciences} \textbf{2020}, 1(1), 7.

\bibitem{Sch12} Scholze, P. Perfectoid Spaces, \textit{Publ. Math. IHES} \textbf{2012} 116, 245–313.  

\bibitem{Sch17} Scholze, P.. Etale cohomology of diamonds, 2017. \url{https://arxiv.org/abs/1709.07343} (accessed January 29, 2021).

\bibitem{FW79} Fontaine, J-M. and Wintenberger, J-P. Extensions alg\'ebrique et corps des normes des extensions APF des corps locaux, \textit{C. R. Acad. Sci. Paris S\'er. A–B} \textbf{1979} 288(8), A441–A444.



\bibitem{bou} ncatlab authors. Bousfield Localization of Model Categories. \url{https://ncatlab.org/nlab/show/Bousfield+localization+of+model+categories} (accessed January 29, 2021).

\bibitem{Cal05} Caldararu, A. \textit{Derived Categories of Sheaves: a Skimming}, 2005. \url{https://arxiv.org/pdf/math/0501094.pdf} (accessed January 29, 2021).

\bibitem{con} ncatlab authors. Condensed Mathematics. Available online: \url{https://ncatlab.org/nlab/show/condensed+mathematics (accessed January 29, 2021)}.


\bibitem{top} ncatlab authors. Topological localization. Available online: \url{https://ncatlab.org/nlab/show/topological+localization (accessed February 2, 2021)}.

\bibitem{site} ncatlab authors. (Infinity,1)-Site. Available online: 
\url{https://ncatlab.org/nlab/show/%28infinity%2C1%29-site (accessed February 2, 2021)}.

\bibitem{reflective} ncatlab authors. Reflective Subcateory. Available online: \url{https://ncatlab.org/nlab/show/reflective+subcategory (accessed February 2, 2021)}.

\bibitem{lur} Lurie, J. \textit{Higher Topos Theory}, Annals of Mathematics Studies 170, Princeton University Press 2009.

\bibitem{Joy08} Joyal, A. Notes on Quasi-Categories, 2008 \url{https://web.math.rochester.edu/people/faculty/doug/otherpapers/Joyal-QC-Notes.pdf} (accessed January 29, 2021).

\bibitem{kan} ncatlab authors. Kan Complex. \url{https://ncatlab.org/nlab/show/Kan+complex} (accessed January 29, 2021).

\bibitem{hom} ncatlab authors. Homotopy Lifting Property. \url{https://ncatlab.org/nlab/show/homotopy+lifting+property} (accessed January 29, 2021).

\bibitem{mod} ncatlab authors. Model category. \url{https://ncatlab.org/nlab/show/model+category} (accessed January 29, 2021).

\bibitem{Gro60} Grothendieck, A., and Dieudonn\'e, J. \textit{\'el\'ements de g\'eom\'etrie alg\'ebrique.} Publications math\'ematiques de l’IHES: Paris, 1960-1967.

\bibitem{tur} Wikipedia authors. Turing degree. \url{https://en.wikipedia.org/wiki/Turing_degree} (accessed January 29, 2021).


\bibitem{juan} Maldacena, J., \textit{The Large N limit of superconformal field theories and supergravity}, Advances in Theoretical and Mathematical Physics. 2 (4): 231–252. arXiv:hep-th/9711200. Bibcode:1998AdTMP...2..231M. doi:10.4310/ATMP.1998.V2.N2.A1.






%


























\end{thebibliography}
\end{document}